\documentclass[11pt]{amsart}     
\usepackage{amssymb,amsthm,graphicx,here,amsmath}
\newtheorem{theorem}{Theorem}

\newtheorem{proposition}{Proposition}
\newtheorem{lemma}{Lemma}

\newtheorem{corollary}{Corollary}
\newtheorem{remark}{Remark}
\newtheorem{assertion}{Assertion}
\newcommand{\ord}{{\rm{ord}}}

\begin{document}
\title{On Alexander polynomials of certain $(2,5)$  torus curves}
\author[M. Kawashima and M. Oka]
{M. Kawashima and M. Oka}
\address{\vtop{
\hbox{Department of Mathematics}
\hbox{Tokyo University of science}
\hbox{wakamiya-cho 26, shinjuku-ku}
\hbox{Tokyo 162-0827}
\hbox{\rm{e-mail}: {\rm M. Kawashima: j1107702@ed.kagu.tus.ac.jp }}
\hbox{\hspace{1.1cm}{\rm M. Oka: oka@rs.kagu.tus.ac.jp }}
}}

\begin{abstract} In this paper, we compute Alexander polynomials of
a torus curve $C$ of  type $(2,5)$, 
 $C:\,f(x,y)=f_2(x,y)^5+f_5(x,y)^2=0$,
 under the assumption that the origin $O$ is the unique inner singularity
and $f_2=0$ is an irreducible conic.
We show  that the Alexander polynomial remains the  same with that of a generic
 torus curve as long as $C$ is irreducible.
\end{abstract}
\maketitle 
\section{Introduction}
A plane curve $C\subset \Bbb P^2$ of degree $pq$ is called a 
$curve \ of \ torus \ type \ (p,q)$ with $p>  q\ge  2$, 
if there is a defining polynomial $F$ of $C$
of the form 
 $F=F_{p}^q+F_{q}^p$,
where $F_p$, $F_q$ are homogeneous polynomials of $X,Y,Z$
 of degree $p$ and $q$ respectively.
A singularity $P\in C$ is called {\em inner} if $F_p(P)=F_q(P)=0$.
Otherwise, $P$ is called an {\em outer singularity.}
A torus curve $C$ is called {\em tame} if it has no
outer singularity.
We assume $O=(0,0)$ hereafter.
In \cite{Kawashima1}, 
the first author  classified the topological
types
of the  germs of inner   singularity of curves of  $(2,5)$ torus type.
In this paper, we are interested in
 the  {\em{Alexander polynomial}} of $C$
which is an important  topological invariant. 
In  the case of irreducible sextics of torus type  $(2,3)$,
there are only 3 possible Alexander polynomials:
$\Delta_{3,2}^j(t)=(t^2-t+1)^{j}$, $j=1,2,3$ (\cite{OkaAtlas}).

A tame torus curve $C$ of type $(p,q)$ is said to be {\em{generic}}
if the associated curves $C_p=\{F_p=0\}$ and $C_q=\{F_q=0\}$ intersect transversely
at $pq$ distinct points.
It is known that  the Alexander polynomial 
of a generic  $C$ is equal to $\Delta_{p,q}(t)$  (\cite{OkaSurvey}) where 
\[
\Delta_{p,q}(t) :=\dfrac{(t^{pq/r}-1)^r(t-1)}{(t^p-1)(t^q-1)},\quad
  r={\text{gcd$(p,q)$}}.
\]
Moreover it is also known that
the Alexander polynomial 
of $C$ is  still equal to $\Delta_{p,q}(t)$,
if  $C$ is tame and $C_p,\,C_q$ intersect at $O$ with intersection multiplicity $pq$ and
$C_p$ is  smooth (\cite{BenoitTu,BenoitTuCorrect}).

Let $C$ be a torus curve of type $(2,5)$
such that 
$C$ has a unique inner singularity, say $O\in C$ (thus
$I(C_2,C_5;O)=10$) and  we assume that $C$ has no outer singularity. 
Then we have shown that there
 are  22 possible singularities for $(C,O)$ under the
assumption that $C_2$ is irreducible (\cite{Kawashima1}).
For 8 classes  among 22 type of singularities,
 $C$ can be  either   irreducible  or  reducible.
We list those 22-singularities below. 
Throughout this paper, 
we use the same notations of singularities
as in  \cite{Kawashima1,okabook}.

\vspace{.3cm}
\noindent
(I) Assume that $C$ is irreducible,  the possibilities are:
\begin{eqnarray*}
&  B_{50,2},\ \ 
   B_{43,2}\circ B_{2,3},\ \
   B_{36,2}\circ B_{4,3},\ \  
   B_{29,2}\circ B_{6,3},\ \ 
   B_{22,2}\circ B_{8,3},\ \
   B_{15,2}\circ B_{10,3},\ \ 
   B_{25,4},\\
&  
   (B_{4,2}^2)^{B_{32,2}+B_{2,2}},\ \ 
   (B_{4,2}^2)^{B_{32,2}+B_{2,2}},\ \
   (B_{6,2}^2)^{B_{23,2}+B_{3,2}},\ \
   (B_{8,2}^2)^{B_{14,2}+B_{4,2}},\  \
   (B_{10,2}^2)^{2B_{5,2}},\\ 
&
   (B_{11,2}^2)^{ B_{6,2}},\ \
   (B_{12,2}^2)^{2{B}_{1,2}},\ \
   (B_{6,2}^2)^{B_{16,2}+B_{1,2}}\circ B_{2,1},\ \
   (B_{8,2}^2)^{B_{7,2}+B_{2,2}}\circ B_{2,1},\\ 
&
   (B_{9,2}^2)^{B_{5,2}}\circ B_{2,1},\ \
    B_{29,2} \circ B_{2,1}  \circ  (B_{2,1}^2)^{B_{k,2}}\ (k=1,2,3,
    5).
\end{eqnarray*}
\noindent
(II) If $C$ is {reducible}, the possibilities are:
\newline
\indent (a) with a line component:
\begin{eqnarray*}
\begin{split}
& \qquad  B_{29,2}\circ B_{6,3},\ \ 
    (B_{6,2}^2)^{B_{16,2}+B_{1,2}}\circ B_{2,1},\ \
    (B_{8,2}^2)^{B_{7,2}+B_{2,2}}\circ B_{2,1},\ \
  (B_{9,2}^2)^{B_{5,2}}\circ B_{2,1},\\ 
 &\qquad   B_{29,2}\circ B_{2,1}\circ (B_{2,1}^2)^{B_{k,2}},\ k=1,2,3,5.
\end{split}
\end{eqnarray*}
\indent (b) with five conics: $B_{20,5}$.

We recall some of the notations.
\begin{eqnarray*}
\begin{split}
B_{p,q}&:\quad x^p+y^q=0,\\
B_{p,q}\circ B_{r,s}&: (x^p+y^q)(x^r+y^s)=0,\ \ q/p< s/r.
\end{split}
\end{eqnarray*}
The singularities listed below
 have  degenerate faces in their Newton boundaries
 and we need one more toric modification for
 their resolutions. See \cite{Kawashima1} for the detail.
\begin{eqnarray*}
& (B_{4,2}^2)^{B_{32,2}+B_{2,2}},\ \ 
   (B_{4,2}^2)^{B_{32,2}+B_{2,2}},\ \
   (B_{6,2}^2)^{B_{23,2}+B_{3,2}},\ \
   (B_{8,2}^2)^{B_{14,2}+B_{4,2}},\  \
   (B_{10,2}^2)^{2B_{5,2}},\\ 
 & (B_{11,2}^2)^{ B_{6,2}},\ \
   (B_{12,2}^2)^{2{B}_{1,2}},\ \
   (B_{6,2}^2)^{B_{16,2}+B_{1,2}}\circ B_{2,1},\ \
   (B_{8,2}^2)^{B_{7,2}+B_{2,2}}\circ B_{2,1},\\ 
&
   (B_{9,2}^2)^{B_{5,2}}\circ B_{2,1},\ \
    B_{29,2} \circ B_{2,1}  \circ  (B_{2,1}^2)^{B_{k,2}}\ (k=1,2,3,
    5).
\end{eqnarray*}

In this paper, we use the method of Libgober \cite{MR85h:14017},
 Loeser-Vaqui\'e
\cite{Loeser-Vaquie} and 
Esnault-Artal (\cite{Artal,Esnault}) for the computation of the 
Alexander polynomials.

\begin{theorem}\label{th-1}
Let $C$ be a tame torus curve of type (2,5).
Suppose that $C$ has a unique inner singularity and 
$C_2$ is irreducible.
Then the  Alexander polynomial $\Delta_{C}(t)$  of $C$ is  given as follows.

{\rm (1)} If $C$ is irreducible (case (I)), then
\[
  \Delta_{C}(t)=\Delta_{5,2}(t)\,\,\,\text{where}\,\,
\Delta_{5,2}(t)=t^4-t^3+t^2-t+1.
\]

{\rm (2)} If $C$ is reducible and have a line component (case (II-a)),
\[
  \Delta_{C}(t)=
    (t-1)(t^4-t^3+t^2-t+1).\]

{\rm (3)}  If $C$ is reducible and $(C,O)\sim  B_{20,5} $ (case (II-b)), 
\[
  \Delta_{C}(t)=
    (t-1)^4(t+1)^4(t^4-t^3+t^2-t+1)^4(t^4+t^3+t^2+t+1)^3.
\]
\end{theorem}

\begin{corollary}
Let $C$ be a tame torus curve of type (2,5) and assume that
there is a degeneration  family $C_t$, $t\in W$ such that 
$C_t\cong C,\,t\ne 0 $ and $C_0$ is an irreducible tame curve with a unique singular
 point $P$ where $W$ is an open neighbourhood of the origin in  $\Bbb C$.
Assume that $(C_0,P)$ is topologically isomorphic to one of 
the above 21 singularities (Case I).
Then the Alexander polynomial $\Delta_C(t)$ is given by $\Delta_{5,2}(t)$.
 \end{corollary}
\begin{corollary}Let $C$ be a tame irreducible torus curve of type (2,5)
 such that  $C_5$ is smooth and $C_2$ is irreducible. 
Then the Alexander polynomial is given by 
$\Delta_{5,2}(t)$.
 \end{corollary}

\section{Alexander polynomial}
Let us consider the affine coordinate
$\Bbb{C}^2=\Bbb{P}^2\setminus \{Z=0\}$  and let $x=X/Z$, $y=Y/Z$.
Let $C$ be a given plane curve of degree $d$ defined by $f(x,y)=0$
and let $O\in C$ be a singular point of $C$ where $O=(0,0)$.
We assume that the line at infinity $\{Z=0\}$ is generic with respect to 
$C$.
 \subsection{Loeser-Vaqui\'e formula}
Consider an embedded  resolution of $(C,O)\subset (\Bbb C^2,O)$, 
$\pi:\tilde{U}\to U$ where $U$ is an open neighborhood of $O$ and
let $E_1,\dots,E_s$ be the exceptional divisors.
Let $(u,v)$ be a local coordinate system centered at $O$
and $k_i$ and $m_i$ be  respective  order of zero of the canonical two form 
$\pi^{*}(du \wedge dv)$ and $\pi^{*}f$ along the divisor $E_i$.
{\em{The adjunction ideal $\mathcal{J}_{O,k,d}$ of $\mathcal{O}_O$}} 
 is defined by 
\[
\mathcal{J}_{O,k,d} 
=\{\phi\in \mathcal{O}_O \ |\  (\pi^{*}\phi)\ge  \sum_{i}([km_i/d]-k_i)E_i\},
 \quad k=1,\dots,d-1
\]
where $[r]$ is the largest integer $n$ such that $n\le r$ for $r \in
  \Bbb{Q}$  (\cite{Artal,Esnault}).

Let $O(j)$ be the set of polynomials in $x,y$ 
whose degree is less than or equal to $j$.
We consider the canonical  mapping 
$\sigma:\,{\Bbb C}[x,y]\to \mathcal{O}_O$ and its restriction:
\[
\sigma_{k}:\,O(k-3)\to \mathcal{O}_O.
\]
Put $V_k(O)=\mathcal{O}_{O}/\mathcal{J}_{O,k,d}$ and
we denote the composition
$O(k-3)\to \mathcal{O}_{O}\to V_k(O)$
by $\bar{\sigma}_{k}$.
Then the Alexander polynomial is given  as follows. 
\begin{lemma}\label{Loeser-Vaquie}
(\cite{MR85h:14017,Loeser-Vaquie,Artal,Esnault})
The reduced Alexander polynomial $\tilde{\Delta}_C(t)$ is given by  the product
\begin{eqnarray}\label{L-V}
\tilde{\Delta}_{C}(t)=\prod_{k=1}^{d-1}\Delta_k(t)^{\ell_k}
\end{eqnarray}
where  $d$ is the degree of $f$,
  $\ell_k$ is the dimension of ${\rm Coker}\, \bar{\sigma}_k$ and
\[\Delta_k(t)=\Big{(}t-\exp(\dfrac{2k\pi i}{d})\Big{)}
\Big{(}t-\exp(-\dfrac{2k\pi i}{d})\Big{)}.
\]
\end{lemma}
We use the method of Esnault-Artal (\cite{Artal})
to compute $\ell_k$.

\begin{remark}
{\rm{The Alexander polynomial $\Delta_C(t)$ is given as 
\[\Delta_C(t)=(t-1)^{r-1}\tilde{\Delta}_{C}(t)\]
where $r$ is the number of irreducible components of $C$ (\cite{OkaSurvey}). }}
Note that for the case of curve of degree $10$.
\begin{eqnarray*}
& \Delta_5(t)=(t+1)^2,\quad \Delta_6(t)\Delta_8(t)=t^4+t^3+t^2+t+1,\ \ 
  \Delta_7(t)\Delta_9(t)=t^4-t^3+t^2-t+1.
\end{eqnarray*}
\end{remark}

\subsection{Pl\"ucker's formula}
We denote  the Milnor number of the singularity of $(C,P)$ by $\mu(C,P)$ and
 the number of locally irreducible components of $(C,P)$ by $r(C,P)$.
We recall the  {\em{generalized Pl\"ucker's formula}}. 
Let $C_1,\cdots, C_r$ be irreducible components
 of $C$ and let $\tilde{C}_1,\dots,\tilde{C}_r$ be their normalizations,   
let $g(\tilde{C}_i)$ be the genus of $\tilde C_i$ and let  $\Sigma(C)$ be 
the singular locus of $C$.
Then
\[
 \chi(\tilde C)=\sum_{i=1}^r (2-2g(\tilde C_i))
=d(3-d)+\sum_{P\in \Sigma(C)}(\mu(C,P)+r(C,P)-1)\le 2r
\]
%
For further  details, 
we refer to \cite{Milnor, Nambabook,MR2107253}.

\section{Outline of the proof of Theorem \ref{th-1}}
We have to consider the following $22$-singularities.
We denote a class of  a singularity  $(C,O)$
which can appear both as an irreducible curve and
 a reducible curve by ${}^{\sharp}(C,O)$. 
In the section 3.2, we will use notation ${}^{irr}(C,O),\,{}^{red}(C,O)$ to
distinguish
 the case of $C$ being  irreducible and 
reducible. 
\begin{eqnarray*}
&  B_{50,2},\  
   B_{43,2}\circ B_{2,3},\ 
   B_{36,2}\circ B_{4,3},\   
   {}^{\sharp}B_{29,2}\circ B_{6,3},\  
   B_{22,2}\circ B_{8,3},\ 
   B_{15,2}\circ B_{10,3},\   
   B_{20,5},\  B_{25,4},\\
&  
   (B_{4,2}^2)^{B_{32,2}+B_{2,2}},\ \ 
   (B_{4,2}^2)^{B_{32,2}+B_{2,2}},\ \
   (B_{6,2}^2)^{B_{23,2}+B_{3,2}},\ \
   (B_{8,2}^2)^{B_{14,2}+B_{4,2}},\  \
   (B_{10,2}^2)^{2B_{5,2}},\\ 
&
   (B_{11,2}^2)^{ B_{6,2}},\ \
   (B_{12,2}^2)^{2{B}_{1,2}},\ \
   {}^{\sharp}(B_{6,2}^2)^{B_{16,2}+B_{1,2}}\circ B_{2,1},\ \
   {}^{\sharp}(B_{8,2}^2)^{B_{7,2}+B_{2,2}}\circ B_{2,1},\\ 
&
   {}^{\sharp}(B_{9,2}^2)^{B_{5,2}}\circ B_{2,1},\ \
   {}^{\sharp} B_{29,2} \circ B_{2,1}  \circ  (B_{2,1}^2)^{B_{k,2}}\ (k=1,2,3,5).
\end{eqnarray*}

\subsection{Divisibility principle and  Sandwich principle}
Suppose we have a degeneration family $C_s,\,s\in  W$ of reducible 
curves such that 
$C_s,\,s\ne 0$  are equisingular family of plane curves. 
Here $W$ is an open neighbourhood of the origin in $\Bbb{C}$. 
We denote this situation as $C_s\overset{s\to 0}\longrightarrow C_0$.
Then we have the divisibility 
  $\Delta_{C_s}(t)\,|\, \Delta_{C_0}(t)$ 
(Theorem 26 of \cite{OkaSurvey}). Suppose that we have two degeneration
series
$C_s\overset{s\to 0}\longrightarrow C_0$ and 
$D_r\overset{r\to 0}\longrightarrow D_0$ such that 
$C_0\cong D_r\,(r\ne 0)$ and assume that 
$\Delta_{C_s}(t)=\Delta_{D_0}(t)$. Then the divisibility implies that 
$ \Delta_{C_s}(t)=\Delta_{C_0}(t)$ (the Sandwich principle).

\subsection{Degeneration series}
Recall that we have the following degeneration series among the
above singularities  (\cite{Kawashima1}): 
\def\mapdown#1{\Big\downarrow\rlap{$\vcenter{\hbox{$#1$}}$}}

\begin{enumerate}
\item Main sequence:
 \[
\begin{matrix}
 B_{50,2}&\longrightarrow & B_{43,2}\circ B_{2,3}&\longrightarrow&
 (B_{4,2}^2)^{B_{32,2}+B_{2,2}}&
\longrightarrow
&
  B_{36,2}\circ B_{4,3}&&\\
&&&&&~~~\searrow {(a)}&\\
& \longrightarrow &   {}^{\sharp}B_{29,2}\circ B_{6,3} 
& \longrightarrow & B_{22,2}\circ B_{8,3}& \longrightarrow &
 B_{15,2}\circ B_{10,3}&\dashrightarrow& B_{20,5}\\
&&\mapdown{(b)}&~~~\searrow (c)&&&\\
\end{matrix} 
\]
where the branched sequences (a)
  from    $(B_{4,2}^2)^{B_{32,2}+B_{2,2}}$ and 
(b), (c) from $B_{29,2}\circ B_{6,3}$ in the main sequence are as follows.

\vspace{.2cm}
\begin{enumerate} 
\item 
      $(B_{4,2}^2)^{B_{32,2}+B_{2,2}}\to
       (B_{6,2}^2)^{B_{23,2}+B_{3,2}}\to 
       (B_{8,2}^2)^{B_{14,2}+B_{4,2}}\to 
       (B_{10,2}^2)^{2B_{5,2}}\to$\newline

      \indent $(B_{11,2}^2)^{ B_{6,2}}\to
       (B_{12,2}^2)^{2{B}_{1,2}}\to
        B_{25,4}$. 

\vspace{.2cm}
\item
  \begin{enumerate}
     \item  ${}^{irr} B_{29,2}\circ B_{6,3}\to 
           {}^{irr} (B_{6,2}^2)^{B_{16,2}+B_{1,2}}\circ B_{2,1}\to
           {}^{irr} (B_{8,2}^2)^{B_{7,2}+B_{2,2}}\circ B_{2,1}\to$
     \vspace{.2cm}\newline\indent
            ${}^{irr}(B_{9,2}^2)^{B_{5,2}}\circ B_{2,1}$.
     \vspace{.2cm}\item ${}^{irr}B_{29,2}\circ B_{6,3}\to  
           {}^{irr} B_{29,2} \circ B_{2,1}\circ (B_{2,1}^2)^{B_{1,2}}\to 
           {}^{irr} B_{29,2} \circ B_{2,1}\circ (B_{2,1}^2)^{B_{2,2}}\to $
     \vspace{.1cm}\newline\indent 
            ${}^{irr}B_{29,2} \circ B_{2,1}\circ (B_{2,1}^2)^{B_{3,2}}\to
           {}^{irr} B_{29,2} \circ B_{2,1}\circ (B_{2,1}^2)^{B_{5,2}}$.
\end{enumerate}

\vspace{.2cm}
\item 
  \begin{enumerate}
     \item ${}^{red}B_{29,2}\circ B_{6,3}\to 
            {}^{red}(B_{6,2}^2)^{B_{16,2}+B_{1,2}}\circ B_{2,1}\to
            {}^{red}(B_{8,2}^2)^{B_{7,2}+B_{2,2}}\circ B_{2,1}\to
            {}^{red}(B_{9,2}^2)^{B_{5,2}}\circ B_{2,1}$.\vspace{.2cm}
     \item 
      ${}^{red}B_{29,2}\circ B_{6,3}\to  
       {}^{red}B_{29,2} \circ B_{2,1}\circ (B_{2,1}^2)^{B_{1,2}}\to 
       {}^{red}B_{29,2} \circ B_{2,1}\circ (B_{2,1}^2)^{B_{2,2}}\to 
       {}^{red}B_{29,2} \circ B_{2,1}\circ (B_{2,1}^2)^{B_{3,2}}\to 
       {}^{red}B_{29,2} \circ B_{2,1}\circ (B_{2,1}^2)^{B_{5,2}}$.\vspace{.2cm}
\end{enumerate}       
\end{enumerate}
\end{enumerate}

The main sequence is obtained through the
 degenerations of  the tangent cone of $C_5$ at $O$,
keeping the irreducibility of $C_2$.
In the last degeneration 
$B_{15,2}\circ B_{10,3}\dashrightarrow B_{20,5}$
of the main sequence,
$C$ degenerates into a reducible curve. 

The branched sequence (a) from $(B_{4,2}^2)^{B_{32,2}+B_{2,2}}$
is obtained  by   degenerating  $(C_5,O)$,
fixing the tangent cone of $C_5$ at $O$.
More precisely,
the tangent cone of $(C_5,O)$ is a line with multiplicity 2
and the generic singularity of $(C_5,O)$ is $A_3$
and the corresponding  degenerations of $(C_5,O)$  are:
\[(C_5,O):\quad B_{4,2}\to B_{6,2}\to B_{8,2}\to B_{10,2}
  \to B_{11,2}\to B_{12,2} \to B_{13,2}.\]
  
\noindent
The branched sequence (b) (respectively, (c))
from ${}^{irr}B_{29,2}\circ B_{6,3}$ (resp. ${}^{red}B_{29,2}\circ B_{6,3}$)
is also obtained by degenerating $(C_5,O)$
fixing the tangent cone of $C_5$ at $O$ (See  \S3.4).

\subsection{Strategy} 
Our strategy is the following.
The singularity $B_{50,2}$ is obtained when $C_2$ and $C_5$
has a maximal contact at $O$ and  $(C_5,O)$ is smooth. 
In this case,  it is known that
$\Delta_C(t)=t^4-t^3+t^2-t+1$ by Theorem 2 of \cite{BenoitTu}.
Hence by virtue of the Sandwich principle, it is enough to show 

(1) the irreducibility of $C$ and 

(2) ${\tilde \Delta}_C(t)=\Delta_{5,2}(t)$ for the case
     $(C,O)$ being one of the following singularities 
which are the end of the degenerations.
\[
B_{15,2}\circ B_{10,3},\ \
B_{25,4},\ \ {}^{\sharp}{(B_{9,2}^2)}^{B_{5,2}}\circ B_{2,1},\ \
{}^{\sharp}B_{29,2} \circ B_{2,1}  \circ  (B_{2,1}^2)^{B_{5,2}}.\] 
By virtue of Lemma \ref{Loeser-Vaquie}, to show ${\tilde
\Delta}_C(t)=\Delta_{5,2}(t)$
is equivalent to show that 

{\em {\rm ($\sharp$)}: $\bar\sigma_{k}:O(k-3)\to V_k(O)$ has one-dimensional cokernel
for $k=7,9$ and surjective for other cases}.
\newline
 So for the proof of 
the assertions
(1) and (2) of Theorem 1, we will actually show the above property $(\sharp)$.

The last singularity $B_{20,5}$ of the main sequence
 appears when $C$ consists of five conics.
We treat this case separately in the later section.

\subsection{Irreducibility of $C$}
Now we will discuss the irreducibility of $C$ using 
 the generalized Pl\"ucker's formula. 
First we show that $C$ is irreducible 
if $(C,O)$ is one of 2 singularities $B_{15,2}\circ B_{10,3}$ and $B_{25,4}$. 

\noindent
{\bf{Case $(C,O)\sim B_{15,2}\circ B_{10,3}$:}}
Note that  the singularities $B_{15,2}$ and $B_{10,3}$
are locally irreducible singularities.
As $\mu(B_{15,2})=14$, $\mu(B_{10,3})=18$
and each singularity appears  for sextics or higher degree curves. 
Thus $C$ must be irreducible, as the degree of $C$ is 10.

\noindent{\bf{Case $(C,O)\sim B_{25,4}$:}}
The singularity $B_{25,4}$ is a  locally
irreducible singularity and thus $C$ is irreducible. 

\vspace{0.2cm}\noindent{\bf{Case $(C,O)\sim B_{29,2}\circ B_{6,3}$:}}
Next  we consider the case $(C,O)\sim B_{29,2}\circ B_{6,3}$ and 
we will show that $C$  can be either irreducible or  reducible.
Recall that the singularity $B_{29,2}\circ B_{6,3}$ 
appears in the case that $C_2$ and $C_5$ satisfies following three
conditions
(\cite{Kawashima1}):

(1) $C_2$ is irreducible and $I(C_2,C_5;O)=10$. 

(2) $(C_5,O)$  has  the multiplicity  3
and  the tangent cone consists of a multiple line $L_1$ of
 the multiplicity 2 and a single line $L_2$.

(3) The conic $C_2$ is tangent to the line  $L_1$ at $O$.

\noindent
Under the condition $I(C_2,C_5;O)=10$, we have generically $(C,O)\sim B_{29,2}\circ B_{6,3}$.
The singularity $B_{29,2}$ is locally irreducible and 
$B_{29,2}$ appears for curves of degree $d\ge 7$ as $\mu(B_{29,2})=28$.
Hence we have four possibilities:
\begin{eqnarray*}
& (1)\ C {\text{: irreducible,}}\quad (2)\  C=D_9\cup D_1,\quad
  (3)\ C=D_8\cup D_2,\quad  (4)\ C=D_7\cup D_3  
\end{eqnarray*}
where $D_d$ is a curve of degree $d$. 
But the cases (3) and (4) are impossible.
Indeed, if $C=D_7\cup D_3$,
then either (a) $(D_7,O)\sim B_{29,2},\,(D_3,O)\sim B_{6,3}$ 
or (b) $(D_7,O)\sim B_{29,2}\circ B_{2,1},\,(D_3,O)\sim B_{4,2}$.
We observe that  $\mu(D_3,O)=10$ in the  case (a) and 
$\mu(D_7,O)=35$ in the case (b) and neither
case is possible by the generalized Pl\"ucker's formula. 
By the same argument, we see that 
the case (3) is impossible.  
Hence we have two possibilities: 

(i) $C$ is irreducible or 

(ii) $C$ consists of  a  line and  a  curve of degree 9. \\
If $C$ has a line component, this line  must be defined by $\{y=0\}$.
%
%
In fact, this case is given by the normal forms of 
$f_2,\,f_5$:
\begin{eqnarray*}
 \begin{split}
f_2(x,y)&=a_{02}\,y^2+(a_{11}\,x+1)\,y-k^2\,x^2,\\
f_5(x,y)&=(t+a_{02}\,b_{04})\,y^5+\phi_4(x)\,y^4+
          \phi_3(x)\,y^3+\phi_2(x)\,y^2+\phi_1(x)\,y-k^5\,x^5
\end{split}
\end{eqnarray*}
where $\phi_1,\phi_2,\phi_3,\phi_4$ take the forms:
\begin{eqnarray*}
\phi_4(x)&=&(a_{02}\,b_{13}-a_{02}^2\,b_{12}+a_{11}\,b_{04})x+b_{04}, \\
\phi_3(x)&=&(b_{13}\,a_{11}-k^2
    \,b_{04}-2b_{12}\,a_{02}\,a_{11}+b_{22}\,a_{02})\,x^2+b_{13}\,x, \\
\phi_2(x)&=&(a_{02}\,k^3+k^2\,a_{02}\,b_{12}-k^2\,b_{13}-b_{12}
           a_{11}^2+b_{22}\,a_{11})\,x^3+b_{22}\,x^2+b_{12}\,x, \\
\phi_1(x)&=&(a_{11}\,k^3+b_{12}\,k^2\,a_{11}-b_{22}\,k^2)\,x^4+(k^3-k^2\,b_{12})\,x^3.
\end{eqnarray*}

\noindent 
The branched sequence (b), (c) in \S 3.2 are obtained by 
degenerating $(C_5,O)$, fixing the tangent cone of $(C_5,O)$
and keeping irreducibility of $C$.

\noindent {\bf{Case $(C,O)\sim B_{20,5}$:}}
This is the last singularity in the main sequence.
We will  show that $C$ can not be irreducible in this case.
As $\mu(B_{20,5})=76$, the number of irreducible components $r$ of $C$
must be at least $5$ by the generalized Pl\"ucker's formula.
On the other hand, 
the singularity $B_{20,5}$ consists of 5 smooth local
components.
Any two components intersects with intersection multiplicity 4.
Thus each local component corresponds to a global component and its
degree must be  2, namely a conic.

\section{Calculation of $\Delta_C(t)$ I: Non-degenerate case}
We divide the calculation of the Alexander polynomial $\Delta_C(t)$
 in two cases,
according to $(C,O)$ being 
non-degenerate or not.  In this section, we treat the first case.
\subsection{Characterization of the adjunction ideal for non-degenerate singularities}
In general,
the computation of the ideal $\mathcal{J}_{O,k,d}$ 
requires an explicit computation of the resolution of the singularity $(C,O)$.
However for the case of non-degenerate singularities,
the ideal $\mathcal{J}_{O,k,d}$ can be obtained combinatorially
by a toric modification.
Let $(u,v)$ be a local coordinate system centered at $O$ such that 
$(C,O)$ is defined by a function germ $f(u,v)$ and
the Newton boundary $\Gamma(f;u,v)$ is non-degenerate.
Let $Q_1,\dots,Q_{s}$ be the primitive weight vectors
which correspond to the faces  $\Delta_1,\dots, \Delta_s$ of $\Gamma(f;u,v)$.
Let $\pi:\tilde{U}\to U$ be the canonical toric modification 
and
let $\hat{E}(Q_i)$ be the exceptional divisor corresponding to $Q_i$.
Recall that the order of zeros of
the canonical two form $\pi^{*}(du\wedge dv)$  along the divisor 
$\hat{E}(Q_i)$ is simply given by $|Q_i|-1$ where 
$|Q_i|=p+q$ for a weight vector $Q_i={^t}{(}p_i,q_i)$ (see \cite{okabook}).
For a function germ $g(u,v)$,
let $m(g,Q_i)$ be the multiplicity of the pull-back 
$(\pi^{*}g)$ on $\hat{E}(Q_i)$.
Then

\begin{lemma}[\cite{OkaAtlas}]\label{OkaCriterion}
A function germ $g\in \mathcal O_O$
 is contained in the ideal $\mathcal{J}_{O,k,d}$ if and only if
$g$ satisfies following condition:
\[
m(g,Q_i)\ge  [\dfrac{k}{d}m(f,Q_i)]-|Q_i|+1,\qquad i=1,\dots,s.
\]  
The ideal $\mathcal{J}_{O,k,d}$ is generated by the monomials
satisfying the above conditions. 
\end{lemma}
We consider the following integers for each singular point $P\in \Sigma(C)$:
\[
 \rho_{k}(P):=\dim V_k(P),\,\,
 \tilde{\rho}(k):=\sum_{P\in \Sigma(C)}\rho_k(P)-\dim{O}(k-3),\,\,
 \iota_{k}(P):=\min_{g\in \mathcal{J}_{P,k,d}}I(g,f;P),
\]
where $V_k(P)=\mathcal{O}_P/\mathcal{J}_{P,k,d}$.
Then the multiplicity $\ell_k$ in the formula  (\ref{L-V}) of Loeser-Vaqui\'e
is given as
\[
\ell_k=\dim {\rm Coker\,} \bar{\sigma}_k
=\tilde{\rho}(k)+\dim {\rm{Ker\,}}\bar{\sigma}_k.
\]
where $\bar\sigma_k$ is defined in \S 2.1.
We consider the integer
 $\sum_{P\in \Sigma(C)}\iota_{k}(P)$.
\begin{proposition} 
If $\sum_{P\in \Sigma(C)}\iota_{k}(P)>d(k-3)$, then 
\begin{enumerate}
\item[(a)] $C$ is irreducible and   $\bar{\sigma}_k$ is injective
and $\ell_k=\tilde{\rho}(k)$ or
\item[(b)] $C$ is reducible.
\end{enumerate}
\end{proposition}
\begin{proof}
Suppose $0\ne g\in {\rm{Ker\,}}\bar{\sigma}_k\subset O(k-3)$.
Then by B\'ezout theorem, we have
\[
d(k-3)\ge  I(G,C)\ge  \sum_{P\in \Sigma(C)}
 I(G,C;P)\ge  \sum_{P\in \Sigma(C)} \iota_k(P)> d(k-3) \]
where $G=\{g=0\}$.
This is an obvious contradiction unless  $g\, |\, f$. 
Thus this implies either $f$ is irreducible and
 $\bar{\sigma}_k$ is injective
 or $f$ is reducible (and $g\,|\,f$).
\end{proof}
\subsection{The singularities $B_{15,2}\circ
B_{10,3}$ and $B_{25,4}$ }
Now
we consider the following  two non-degenerate singularities
 $B_{15,2}\circ
B_{10,3}$, and $B_{25,4}$
which appear as the last singularities of the respective degenerations
with $C$ being irreducible. 
We assume that we have chosen local
analytic coordinates $(u,v)$ so that
\begin{equation*}
\begin{split}
B_{15,2}\circ B_{10,3}&: \,\, f(u,v)=u^{25}+u^{10}v^2+v^5+\text{(higher terms)}, \\
\qquad B_{25,4}&:\,\,f(u,v)=u^{25}+v^4+\text{(higher terms)}.
\end{split}
\end{equation*}
%
The local data are given by the following tables.
\begin{center}
 $B_{15,2}\circ B_{10,3}$\ :\hspace{1cm}
{\renewcommand\arraystretch{1.3}
\begin{tabular}{|c|c|c|c|}
\hline
$k$\ & ${\mathcal{J}}_{O,k,10}$ & $\rho_{k}(O)$ & $\iota_{k}(O)$ \\
\hline
$3$\ & $\langle u,v \rangle$ &1 & 5 \\                            
\hline
$4$\ & $\langle u^3,v \rangle$ &3 &15 \\                            
\hline
$5$\ & $\langle u^5,uv, v^2 \rangle$ &6 &23 \\                            
\hline
$6$\ & $\langle u^{7},u^3v,v^2 \rangle$ &10 &33 \\                            
\hline
$7$\ & $\langle u^{10},u^5v,uv^2,v^3 \rangle$ &16 &43 \\                            
\hline
$8$\ & $\langle u^{12},u^6v,u^3v^2,v^3 \rangle$ &21 &52 \\
\hline
\ $9$\ & \ $\langle u^{15},u^{8}v,u^5v^2,uv^3,v^4 \rangle$ \ &29 &63   \\
       \hline
\end{tabular}
}

\vspace{0.5cm}
$B_{25,4}$\ :\qquad \hspace{1.5cm}
{\renewcommand\arraystretch{1.3}
\begin{tabular}{|c|c|c|c|}
\hline
$k$\ & ${\mathcal{J}}_{O,k,10}$ & $\rho_{k}(O)$ & $\iota_{k}(O)$ \\
\hline
$3$\ & $\langle u,v \rangle$ &1 & 4 \\                            
\hline
$4$\ & $\langle u^3,v \rangle$ &3 &12 \\                            
\hline
$5$\ & $\langle u^6,v \rangle$ &6 &24 \\                            
\hline
$6$\ & $\langle u^{8},u^2v,v^2 \rangle$ &10 &32 \\                            
\hline
$7$\ & $\langle u^{11},u^5v,v^2 \rangle$ &16 &44 \\                            
\hline
$8$\ & $\langle u^{13},u^7v,uv^2,v^3 \rangle$ &21 &52 \\                            
\hline
\ $9$\ & \ $\langle u^{16},u^{10}v,u^3v^2,v^3 \rangle$ \ &29 &62 \\                 \hline
\end{tabular}
}
\end{center}

\noindent
{\bf Case} $(C,O)\sim B_{15,2}\circ B_{10,3}$ and $B_{25,4}$.
In this case, we have the inequalities $\iota_{k}(O)>10(k-3)$ for all
 $k=3,\dots, 9$ 
by the local data.
Hence $\bar{\sigma}_k$ is injective for all $k$ by Proposition 1 and we 
obtain the property ($\sharp$):
\begin{equation*}
\ell_k=\tilde{\rho}(k)=
  \begin{cases}
    1&  k=7,9,\\
    0&  k\ne 7,9.
  \end{cases}
\end{equation*}
Therefore $\Delta_{C}(t)=\Delta_{5,2}(t)=t^4-t^3+t^2-t+1$.

\subsection{ Exceptional case:{\bf  $(C,O)\sim B_{20,5}$}}
In this section, we consider the last singularity $B_{20,5}$
which takes place for reducible $C$.
Recall that $C$ is a union of five conics. 
We assume that we have chosen local coordinates $(u,v)$ so that 
$(C,O)$ is defined by 
\[
 B_{20,5}:\ f(u,v)=u^{20}+v^5+\text{(higher terms)},
\]   
where we ignore the coefficients of the monomials and other monomials
 corresponding 
to other integral points on the Newton boundary.

{\renewcommand\arraystretch{1.3}
\vspace{0.5cm}
  $B_{20,5}$\ :\hspace{2.1cm}
\begin{tabular}{|c|c|c|c|}
\hline
$k$\ & ${\mathcal{J}}_{O,k,10}$ & $\rho_{k}(O)$ & $\iota_{k}(O)$ \\
\hline
$3$\ & $\langle u^2,v \rangle$ & 2 & 10 \\                            
\hline
$4$\ & $\langle u^4,v \rangle$ &4 &20 \\                            
\hline
$5$\ & $\langle u^6,u^2v, v^2 \rangle$ &8 &30 \\                            
\hline
$6$\ & $\langle u^{8},u^4v,v^2 \rangle$ &12 &40\\                            
\hline
$7$\ & $\langle u^{10},u^6v,u^2v^2,v^3 \rangle$ &18 &50 \\                            
\hline
$8$\ & $\langle u^{12},u^8v,u^4v^2,v^3 \rangle$ &24 &60 \\                            
\hline
\ $9$\ & \ $\langle u^{14},u^{10}v,u^6v^2,u^2v^3,v^4 \rangle$ \ &32 &70 \\                            
\hline
\end{tabular}
}

\vspace{0.5cm}
\noindent
Again 
we have the inequalities $\iota_{k}(O)-10(k-3)>0$ for all
 $k=3,\dots, 9$.
We claim that {\em  $\bar{\sigma}_k$ is injective for all $k$}.
In fact,  assuming   $0\ne g\in{\rm{Ker\,}} \bar{\sigma}_k$, 
we have $g\,|\,f$ by  the proof of Proposition 1
and 
this means $g$ is a union of conics 
which are components of $f$. 
Consider the factorization $f=h_1h_2h_3h_4h_5$ 
where $\{h_i=0\}$ is a 
smooth conic component of $C$.
Then we may assume that 
\begin{eqnarray*}
\begin{split}
f &\overset{\sigma}\longmapsto u^{20}+v^5+\text{(higher terms)},\quad 
h_i \overset{\sigma}\longmapsto u^4+\zeta^i v+\text{(higher terms)},\ 
i=1,\dots,5
\end{split}
\end{eqnarray*}
where $\zeta=\exp(\pi i /5)$.
Thus suppose that  $g=h_{i_1}\cdots h_{i_j}$. Then $2j\le k-3$ or
$j\le [\frac{k-3}2]$ and 
$\sigma_k(g)$ must contain $v^j$ with a non-zero coefficient.
This implies that
$j\le 0,0,1,1,2,2,3$ for $k=3,4,\dots,9 $ respectively.
On the other hand, $v^j\in \mathcal{J}_{O,k,10}$
implies from the table of $B_{20,5}$ that
$ j\ge 1,1,2,2,3,3,4$ for $ k=3,\dots,9$ respectively.
This gives an obvious contradiction.
Hence  we have 
\begin{equation*}
\ell_k=\tilde{\rho}(k)=
  \begin{cases}
    1&  k=3,4,\\
    2&  k=5,6,\\
    3&  k=7,8,\\
    4&  k=9.
  \end{cases}
\end{equation*}
Therefore by the formula (\ref{L-V}) in Lemma \ref{Loeser-Vaquie}  we obtain the equality:
$$\Delta_{C}(t)=(t-1)^4(t+1)^4(t^4-t^3+t^2-t+1)^4(t^4+t^3+t^2+t+1)^3.$$

\section{Calculation of $\Delta_C(t)$, II: Degenerate cases }
Next we calculate the Alexander polynomial of following 
two degenerate singularities:
\begin{itemize}
\item ${(B_{9,2}^2)}^{B_{5,2}}\circ B_{2,1}$: this is the last
      singularity of the sequence of (b-i) or (c-i). 
\item $ B_{29,2} \circ B_{2,1}  \circ (B_{2,1}^2)^{B_{5,2}}$:
this is the last singularity of the sequence of (b-ii) or (c-ii). 
\end{itemize}
\subsection{Characterization of the adjunction ideal for degenerate cases}
For degenerate singularities,
we proceed several toric modifications to obtain their resolutions.
Consider an embedded  resolution of $(C,O)\subset (\Bbb C^2,O)$, 
$\pi:\tilde{U}\to U$ where $U$ is an open neighborhood of $O$ and
let $E_1,\dots,E_s$ be the exceptional divisors.
 We put the ideal
$\bar{\mathcal{J}}_{O,k,d}$ of $\mathcal{O}_O$  
\[
\bar{\mathcal{J}}_{O,k,d}:= 
\langle M\in \mathcal{O}_O \ | \ 
       {\text{$M$ : monomial}},\ (\pi^{*}M)\ge  \sum_{i}([km_i/d]-k_i)E_i\ \rangle,\quad 1\le k \le d-1.
\]
In general,
$
\bar{\mathcal{J}}_{O,k,d}\subset \mathcal{J}_{O,k,d}$ and
$\bar{\mathcal{J}}_{O,k,d}= \mathcal{J}_{O,k,d}$
if $(C,O)$ is non-degenerate from Lemma 2.
If $(C,O)$ is degenerate singularity,
there exist several other (non-monomial)  polynomials
 $h_i,\,i=1,\dots,r$ such that 
$h_i\in  \mathcal{J}_{O,k,d}\setminus \bar{\mathcal{J}}_{O,k,d}$ and 
\[
 \mathcal{J}_{O,k,d}=\langle  M ,h_i \,|\,
  M\in\bar{\mathcal{J}}_{O,k,d},\,
i=1,\dots,r  \rangle.
\]

\subsubsection{\bf Formulation of the multiplicities.}
We recall how the  multiplicities of the pull-back of  a function after
toric modifications along the exceptional divisors can be computed. 

Let $D=\{g=0\}$ be a plane curve and let $P\in D$ be a singular point.
Suppose that its Newton boundary $\Gamma(g;u,v)$ consists of $m$-faces
$\Delta_1,\dots, \Delta_m$
where $(u,v)$ is a local coordinates centered at $P$.
Then the face function of $g$  with respect to a face $\Delta_i$ 
 takes the form:
\[g_{\Delta_i}(u,v)=c\,u^{w_i}v^{t_i}\prod_{j=1}^{k_i}
           (v^{a_i}-\gamma_{i,j}u^{b_i})^{\nu_{i,j}},\quad c\ne 0\]
where $P_i={}^t(a_i,b_i)$ is the weight vector 
corresponding to $\Delta_i$.
Let $\{E_0,P_1,\dots,P_m, E_{2}\}$
 be the vertices of the dual Newton diagram $\Gamma^{*}(g;u,v)$ where
 $E_1={}^t(1,0)$ and $E_2={}^t(0,1)$.
Let $\pi_1:X_1\to \Bbb{C}^2$ be the toric modification associated with 
$\{\Sigma^{*}_1,(u,v),P\}$ where 
$\Sigma^{*}_1=\{E_1,Q_1,\dots, Q_{m'},E_2 \}$ is the canonical regular simplicial cone subdivision
of $\{E_1,P_1,\dots,P_m, E_{2}\}$ (\cite{okabook}).  
Then we can write the divisor   $(\pi_1^{*}g)$ as 
\[
(\pi_1^{*}g)=\tilde{D}
          +\sum_{s=1}^{m'}m(g,Q_s)\hat{E}(Q_s)
\]
where $\tilde{D}$ is the strict transform of $D$ and
$\hat{E}(Q_j)$ is the exceptional divisor
corresponding to the vertex $Q_j$. 
We assume that $P_i=Q_{\nu_i}$ for $i=1,\dots, m$.
Then the exceptional divisors
$\hat E(Q_{\nu_i})=\hat E(P_i)$ intersects with
the strict transform $\tilde D$.
We take the toric coordinates 
$(\Bbb C_{\sigma_{\nu_i}}^2,(u_i,v_i))$ 
where $\sigma_{\nu_i}={\rm{Cone}}\,(Q_{\nu_i},Q_{\nu_i+1})$ so that
$\{u_i=0\}$ defines $\hat{E}(Q_{\nu_i})\cap \Bbb C_{\sigma_{\nu_i}}^2$.
Then $\tilde D$ and the total transform $\pi_1^* D$ are 
 defined in this coordinate as
\begin{eqnarray*}
\begin{split}
&\tilde D:\quad 
\tilde{g}(u_i,v_i)= c_i\,(v_i-\gamma_{i,j})^{\nu_{i,j}}+R(u_i,v_i)=0,\quad c_i\ne 0\\
& \pi_1^* D:\quad
 \pi^*g(u_i,v_i)=u_i^{d(P_i;g)}v_i^{d(Q_{\nu_{i}+1};g)}\,\tilde{g}(u_i,v_i)
\end{split}
\end{eqnarray*}
where $R\equiv 0$ modulo $(u_i)$.  
Thus $\xi_{i,j}:=(0,\gamma_{i,j})$
is the intersection points of
$\tilde{D}$ and  $\hat{E}(Q_{\nu_i})$ for $j=1,\dots,k_i$.
 We take an admissible translated  coordinates $(u_i,v_i')$ with 
$v_i'=v_i-\gamma_{i,j}+h(u_i)$ 
in an open  neighbourhood of $\xi_{i,j}$ where $h$ is a suitable
polynomial
with $h(0)=0$.
Suppose that $(\tilde{D},\xi_{i,j})$ has a 
non-degenerate singularity with respect to the coordinates $(u_i,v_i')$
and suppose that the Newton boundary has a unique
 face $\Delta_{i,j}$ for $j=1,\dots,k_i$.
(For our purpose, this case is enough to be considered.) 
Let $S_{i,j}={}^t(s_{i,j},t_{i,j})$
be the primitive dual vector 
which corresponds to the face $\Delta_{i,j}$
and assume the germ  $(\tilde{D},\xi_{i,j})$ is equivalent to
the Brieskorn singularity $B_{c_{i,j},d_{i,j}}$ with 
$t_{i,j}c_{i,j}=s_{i,j}d_{i,j}$.
This means the dual Newton diagram $\Gamma^*(\tilde{g};u_i,v_i')$ is
given by $\{E_1,S_{i,j},E_2\}$. 
 
We  take the canonical
regular subdivision 
$\Sigma_{i,j}^*$ of 
$\Gamma^*(\tilde{g};u_i,v_i') $. Put 
\[
 \Sigma_{i,j}^{*}=\{T_{i,j,0},T_{i,j,1},\dots,
 T_{i,j,m_j},T_{i,j,m_j+1}\},\,
T_{i,j,0}=E_1,T_{i,j,m_j+1}=E_2.
\]
We may assume $S_{i,j}=T_{i,j,k_0}$ for some $k_0\in \{1,\dots, m_j\}$.
At each point $\xi_{i,j}$,
 we take the  toric modification $\pi_{ij}: X_{ij}\to X_1$
with respect to $\{\Sigma_{i,j}^{*},(u_i,v_i'),\xi_{i,j}\}$.
These modifications are compatible each other and  let
$\pi_2: X_2\to X_1$ be the composition of these modifications for every
$i,j$ so that the exceptional divisors of $\pi_2$ are
bijectively corresponding to the vertices of
$\Sigma_{i,j}^{*},\,i=1,\dots,m,\, j=1,\dots, k_i$. 
What is necessary to be
checked are the multiplicities of $\pi^* g$ and 
$\pi^*(du\wedge dv)$ along the exceptional divisors
$\hat{E}( T_{i,j,k})$ where  $\pi: X_2\to \Bbb C^2$ is the composition of 
$\pi_2:X_2\to X_1$ and $\pi_1:X_1\to \Bbb C^2$.
Then we can write:
\begin{equation*}
\begin{split}    
(\pi^*g)&=\tilde{D}+\sum_{s=1}^{m'}m(g,Q_s)\hat{E}(Q_s)+
\sum_{i=1}^{m}
\sum_{j=1}^{k_i}\sum_{k=1}^{m_j}m(g,T_{i,j,k})\hat{E}(T_{i,j,k})
.\\
(\pi^*K)&=\sum_{s=1}^{m'}k(Q_s)\hat{E}(Q_s)+
\sum_{i=1}^{m}
\sum_{j=1}^{k_i}\sum_{k=1}^{m_j}k(T_{i,j,k})\hat{E}(T_{i,j,k})
\end{split}
\end{equation*}
where $K=du\wedge dv$ is the canonical two form in the base space.
\begin{lemma}
Under the above situations,
the multiplicities  are given as follows.
Put $T_{i,j,k}={}^t(\varepsilon_{i,j,k},\eta_{i,j,k})$.
\begin{enumerate}
\item The  multiplicities $m(g,P_i)$, $m(g,T_{i,j,k})$ of $\pi^{*}g$ 
along the divisors $\hat{E}(P_i)$ and  
$\hat{E}(T_{i,j,k})$ are given by 
\begin{eqnarray*}
\begin{split}
m(g,P_i)&=d(P_i,g),\,\,
m(g,T_{i,j,k})&=\varepsilon_{i,j,k}m(g,P_i)+d(T_{i,j,k},\tilde{g}). 
\end{split}
\end{eqnarray*}
\item
 The  multiplicities $k(Q_s)$, $k(T_{i,j,k})$ of 
the pull-back of the canonical two form $K=du\wedge dv$ along the  divisors  $\hat{E}(Q_s)$ and  
$\hat{E}(T_{i,j,k})$ are given by 
\[
k(Q_s) =|Q_s|-1,\quad  
k(T_{i,j,k})=|T_{i,j,k}|-1+\varepsilon_{i,j,k}k(P_i)
\]
where $|{}^t(a,b)|=a+b$.
\end{enumerate}
\end{lemma} 
The proof follows easily from Theorem 3.8 and
Proposition 7.2, Chapter III of \cite{okabook}.
\subsection{Generalization of Lemma \ref{OkaCriterion}}
\begin{lemma}
Under the above assumptions,
a  germ $\varphi \in  \mathcal{O}_P$ is contained in the ideal $\mathcal{J}_{P,k,d}$
 if and only if $\varphi$ satisfies:
\begin{enumerate}
\item $m(\varphi,P_{i})\ge [\dfrac{k}{d}m(g,P_i)]-k(P_i)$ 
                                   for  $i=1,\dots,m$, and
\item $m(\varphi,S_{i,j})\ge [\dfrac{k}{d}m(g,S_{i,j})]-k(S_{i,j})$
                                   for $j=1,\dots, k_i$.
\end{enumerate}
\end{lemma}
Note that there are no conditions on other exceptional divisors 
$\hat{E}(T_{i,j,k})$. 
\begin{proof}
The proof is almost parallel  to that of Lemma 2 of \cite{OkaAtlas}.
Assume that $\varphi$ satisfies the conditions (1) and (2).
It is enough to show that 
\[
(2-bis)\ \  m(\varphi,T_{i,j,k})\ge [\dfrac{k}{d}m(g,T_{i,j,k})]-k(T_{i,j,k}),\quad
  j=1,\dots,k_i,\ k=1,\dots,m_j. 
\] 
Note that the condition (2) is equivalent to
\begin{enumerate}
\item[(2)'] $m(\varphi,S_{i,j})> 
\dfrac{k}{d}m(g,S_{i,j})-(|S_{i,j}|+s_{i,j}k(P_i))$
                                   for $j=1,\dots, k_i$.
\end{enumerate}
 First we observe that 
$m(g,T_{i,j,0})=m(g,P_i)$ and $m(g,T_{i,j,m_j+1})=0$. 
Take $T_{i,j,k}$ for $k< k_0$ for example. We can write 
$T_{i,j,k}=\alpha_k S_{i,j}+ \beta_k T_{i,j,0}$ 
 for some positive rational numbers $\alpha_k,\, \beta_k$. Note that 
\begin{equation*}
\begin{split}
|T_{i,j,k}|&=\alpha_k |S_{i,j}|+ \beta_k |T_{i,j,0}|=\alpha_k |S_{i,j}|+\beta_k,\\
m(g,T_{i,j,k})&=\alpha_k m(g,S_{i,j})+ \beta_k m(g, T_{i,j,0}),\\
\end{split}
\end{equation*}
Here the second equality follows as
$\Delta(S_{i,j},\pi_1^*g)\cap\Delta(T_{i,j,0},\pi_1^*g)\ne \emptyset$
 by the admissibility of the canonical subdivision $\Sigma_{i,j}^*$.
Thus we have
\begin{equation*}
\begin{split}
m(\varphi,T_{i,j,k})
& \ge \alpha_k m(\varphi,S_{i,j})+\beta_k m(\varphi,T_{i,j,0})\\
& >\alpha_k\left(\dfrac{k}{d}m(g,S_{i,j})-(|S_{i,j}|+s_{i,j}k(P_i))\right)+
\beta_k\left(\dfrac{k}{d}m(g,T_{i,j,0})-(1+ k(P_i))
\right)\\
& =\dfrac{k}{d}m(g,T_{i,j,k})-(|T_{i,j,k}|+\varepsilon_{i,j,k}k(P_i))
 \end{split}
\end{equation*} 
as $\varepsilon_{i,j,k}=\alpha_ks_{i,j}+\beta_k$ by the equality 
 $T_{i,j,k}=\alpha_k S_{i,j}+ \beta_k T_{i,j,0}$.
This inequality is equivalent:
\[
m(\varphi,T_{i,j,k})\ge [\dfrac{k}{d}m(g,T_{i,j,k})]-k(T_{i,j,k}).
\]
For $T_{i,j,k}$ with $k> k_0$, the argument is similar.
Hence we have $\varphi \in \mathcal{J}_{P,k,d}$.
\end{proof}

Now we consider the ideal $\mathcal{J}_{P,k,d}$ 
in more detail.
Take $\varphi \in \mathcal{O}_P$.
We compute the multiplicity of $\varphi$
along the divisors $\hat{E}(P_i)$ and   
$\hat{E}(S_{i,j})$.
We divide our consideration into the two cases:
\begin{enumerate}
\item $\varphi$ is a monomial, 

\item $\varphi$ is a polynomial (non-monomial).  
\end{enumerate}
First we see the case (1) and
we put $\varphi(u,v)=u^{\alpha}v^{\beta}$. 
As $\pi_1^* \varphi$ is also a monomial in
$u_i,v_i$, 
 we can check easily following 
\[
m(\varphi,P_i)=d(P_i,\varphi)=a_i\alpha+b_i\beta, \quad 
m(\varphi,S_{i,j})=s_{i,j}m(\varphi,P_i).
\]  

Next we consider the case (2). 
We can write $\varphi(u,v)=\varphi_{P_i}(u,v)+R(u,v)$
where $R(u,v)$ consist of monomials of degree strictly
greater than $d(P_i,\varphi)$.
If $\Delta(\varphi,P_i)$ is zero dimensional,
then the multiplicities $m(\varphi,P_i)$ and $m(\varphi,S_{i,j})$ are equal to
that of the monomial $\varphi_{P_i}(u,v)$.
If $\Delta(\varphi,P_i)$ is one dimensional,
then the face function $\varphi_{P_i}(u,v)$ can be written by 
\[
\varphi_{P_i}(u,v)=c_i\,u^{\alpha}v^{\beta}\prod_{j=1}^{\kappa_i}
(v^{a_i}-\delta_{i,j}u^{b_i})^{\mu_{i,j}},\quad c_i,\ \delta_{i,j}\ne0.
\]
Then the multiplicities $m(\varphi,P_i)$ 
is given by 
\begin{eqnarray*}
 m(\varphi,P_i)&=&a_i\alpha+b_i\beta+a_ib_i\,\sum_{j=1}^{\kappa_i}\mu_{i,j}.
\end{eqnarray*}
In the admissible translated  coordinates $(u_i,v_i')$,
the function $\pi_1^*\varphi$ is written by
\begin{eqnarray*}
 \begin{split}
 \pi_1^*\varphi(u_i,v'_i)&=c_iu_i^{m(\varphi,P_i)}
                           \tilde{\varphi}(u_i,v'_i),\\
 \tilde{\varphi}(u_i,v'_i)&=\prod_{j=1}^{\kappa_i}
  (v'_i+(\gamma_{i,j}-\delta_{i,j})-h(u_i))^{\mu_{i,j}}
  +\tilde{R}(u_i,v'_i)
 \end{split}
\end{eqnarray*}
  where $\tilde{R}(u_i,v'_i)\equiv 0 \mod (u_i)$.
Thus we obtain
\begin{eqnarray*}
 m(\varphi,S_{i,j})&=&
\begin{cases}
s_{i,j}m(\varphi,P_i)& \text{if $\delta_{i,j}\ne \gamma_{i,j}$ for all $j$},\\
s_{i,j}m(\varphi,P_i)+d(S_{i,j},\tilde{\varphi})&
                       \text{if $\delta_{i,j}= \gamma_{i,j}$ for some $j$.}
\end{cases}
\end{eqnarray*}
Note that 
the multiplicity $d(S_{i,j},\tilde{\varphi})$ depends on
the form $h$, $R$ and $S_{i,j}$. 

\subsection{The case of $(B_{9,2}^2)^{B_{5,2}}\circ B_{2,1}$}
By the local classification in  \cite{Kawashima1},
this singularity $(B_{9,2}^2)^{B_{5,2}}\circ B_{2,1}$ 
appears when the associated curves
$C_2$ and $C_5$ satisfies following conditions:
\begin{enumerate}
\item $C_2$ is irreducible and $I(C_2,C_5;O)=10$.

\item The multiplicity of $(C_5,O)$ is 3
and the  tangent cone of $C_5$ consists of a line $L_1$ with multiplicity 2 and
             a single line $L_2$.

\item The conic $C_2$ is tangent to the line  $L_1$  at $O$.
\end{enumerate}  
Suppose that $C_2$ and $C_5$ satisfies the above conditions.
Then we may assume that the defining polynomials of  $C_2$ and $C_5$ 
are the following forms:
\begin{eqnarray*}
\begin{split}
f_2(x,y)&=y+a_{20}x^2+a_{11}xy+a_{02}y^2,\ a_{20}\ne 0,\\
f_5(x,y)&=b_{05}y^5+((a_{02}^2b_{12}+a_{11}b_{04})x+b_{04})y^4
          +((2b_{12}a_{02}a_{11}+a_{20}b_{04})x^2+2a_{02}b_{12}x)y^3\\
         & \qquad \qquad \qquad+((2a_{20}a_{02}b_{12}+b_{12}a_{11}^2)x^3+2b_{12}a_{11}x^2+b_{12}x)y^2\\
         &  \qquad \qquad \qquad \qquad+(2a_{11}b_{12}a_{20}x^4+2b_{12}a_{20}x^3)y+a_{20}^2b_{12}x^5
\end{split}
\end{eqnarray*}
where $b_{12}\ne 0$ and  $a_{20}+b_{12}^2\ne 0$ in general. 
If $a_{20}+b_{12}^2=0$, 
$(C,O)$ has the same type of singularity but $C$ is not irreducible and has a line component which is defined by $\{y=0\}$.
Now we take a local coordinates $(u,v)$  of the following type so that
\begin{eqnarray*}
\begin{split}
 x=u,&\quad y=v+\varphi(u),\quad \varphi(u)=-a_{20}u^2+\cdots,\ a_{20}\ne 0,\\
f_2(u,v+\varphi(u))&=v+c_1u^5+\cdots, \\
f_5(u,v+\varphi(u))&=b_{04}v^4+b_{12}uv^2+c_2u^{10}+{\text{(higher
 terms)}}, 
\, b_{12}, c_2\ne 0, \\
 f(u,v+\varphi(u))&=v^5+u^2(b_{12}v^2+c_2u^9)^2+{\text{(higher terms)}}. 
\end{split}
\end{eqnarray*}
Then the Newton boundary  $\Gamma(f;u,v)$ consists of two faces 
$\Delta_i\ (i=1,2)$ so that 
the respective  face functions are given by
\[
f_{\Delta_1}(u,v)=v^4(v+b_{12}^2u^2),\quad 
f_{\Delta_2}(u,v)=u^2(b_{12}v^2+c_2u^9)^2.
\] 
Note that  $f(u,v)$ is degenerate on $\Delta_2$.
We take the  canonical toric modification  $\pi_1:X_1\to \Bbb{C}^2$
with respect to 
$\{\Sigma_1^{*},(u,v),O \}$ where
$\Sigma_1^{*}$ is the canonical regular simplicial cone subdivision with
vertices $\{E_1,Q_1,\dots, Q_6,E_2\}$ where
\begin{equation*}
Q_1 = \begin{pmatrix}
         1  \\
         1 
         \end{pmatrix},\
Q_2 = \begin{pmatrix}
         1  \\
         2 
         \end{pmatrix},\
Q_3 = \begin{pmatrix}
         1  \\
         3 
         \end{pmatrix},  \
Q_4 = \begin{pmatrix}
         1  \\
         4 
         \end{pmatrix},  \
Q_5 = \begin{pmatrix}
         2  \\
         9 
         \end{pmatrix},  \
Q_6 = \begin{pmatrix}
         1  \\
         5 
         \end{pmatrix}                     
\end{equation*}
and the weight vectors $Q_2$ and $Q_5$  
correspond to the faces $\Delta_1$ and $\Delta_2$ respectively.
Then we can write the divisor $(\pi_1^{*}f)$ as
\[(\pi_1^{*}f)= \tilde{C}+\sum_{i=1}^6m(f,Q_i)\hat{E}(Q_i),\]
where $\tilde{C}$ is the strict transform of $C$
and  intersects only with the exceptional divisors 
$\hat E(Q_2)$ and $\hat E(Q_5)$.
We can see that $\tilde{C}$ is smooth and intersects transversely at 
$\tilde{C} \cap \hat{E}(Q_2)$ but 
$\tilde{C}$ has the singularity at  the intersection 
$\tilde{C} \cap \hat{E}(Q_5)$. Put 
$\xi= \tilde{C} \cap \hat{E}(Q_5)$. 
In the toric coordinates $(u_1,v_1)$ of
$\Bbb{C}_{\tau}^2$ with $\tau={\rm{Cone}}\,(Q_5,Q_6)$
(see \cite{okabook} for the notations), $\xi=(0,-c_2/b_{12})$.
To see the singularity $(\tilde{C},\xi)$,
we take the admissible translated toric coordinates $(u_1,v_1')$ with
$v_1'=v_1+c_2/b_{12}+h(u_1)$ where $h$ take the form $h(u_1)=q_1u_1+q_2u_1^2$.
Then we can see that 
$\pi_1^*f(u_1,v_1')=c\,u_1^{40}({v_1'}^2+\beta
u_1^5+{\text{(higher terms)}})$ and
  $(\tilde{C},\xi)\sim B_{5,2}$.
Now we take the second  toric modification $\pi_2:X_2\to X_1$
with respect to $\{\Sigma_2^{*},(u_1,v'_1),\xi\}$ where
$\Sigma_2^{*}$ is the canonical regular 
simplicial cone subdivision with vertices
$\{E_1,T_1,\dots, T_4,E_2\}$ where
\begin{equation*}
T_1 = \begin{pmatrix}
         1  \\
         1 
         \end{pmatrix},\
T_2 = \begin{pmatrix}
         1  \\
         2 
         \end{pmatrix},\
T_3 = \begin{pmatrix}
         2  \\
         5 
         \end{pmatrix},  \
T_4 = \begin{pmatrix}
         1  \\
         3
         \end{pmatrix},  
\end{equation*}
and the weight vector $T_3$ 
corresponds to the unique face of $\pi_1^{*}f(u_1,v_1')$.
Note also the exceptional divisor which corresponds to $E_1$ is 
nothing but the exceptional divisor $\hat{E}(Q_5)$ in the previous 
modification $\pi_1$.
Then we have 
\begin{multline*}
(\pi^{*}f)=5\hat{E}(Q_1)+10\hat{E}(Q_2)+14\hat{E}(Q_3)+18\hat{E}(Q_4)
           +40\hat{E}(Q_5)+20\hat{E}(Q_6)
       \\ +42\hat{E}(T_1)+44\hat{E}(T_2)+90\hat{E}(T_3)+45\hat{E}(T_4)
\end{multline*}
\begin{multline*}
(\pi^{*}K)=\hat{E}(Q_1)+2\hat{E}(Q_2)+3\hat{E}(Q_3)+4\hat{E}(Q_4)
           +10\hat{E}(Q_5)+5\hat{E}(Q_6)\\ 
           +11\hat{E}(T_1)+12\hat{E}(T_2)+26\hat{E}(T_3)+13\hat{E}(T_4)
\end{multline*}
and we consider two polynomials
$h_2(u,v)$ and $r_2(u,v)$ which are defined by 
$h_2(u,v)=b_{12}v^2+c_2u^9$ and
$r_2(u,v)=h_2(u,v)-\frac{q_1b_{12}^2}{c_2}u^5v$
  Then we can see by a direct computation
  \begin{eqnarray*}
& \pi_1^*h_2(u_1,v_1')=u^{18}(d_3v'_1+d_4u_1+\text{(higher terms)}),\\   
& \pi_1^*r_2(u_1,v_1')=u^{18}(d'_3v'_1+d'_4u_1^2+\text{(higher terms)}),\\  
&m(h_2,Q_2)=4,\ \  m(h_2,Q_5)=18,\ \  m(h_2,T_3)=38,\\
& m(r_2,Q_2)=4,\ \  m(r_2,Q_5)=18,\ \  m(h_2,T_3)=40.
  \end{eqnarray*}
\begin{assertion}
The adjunction ideals $\mathcal{J}_{O,k,10}$ are given by 
\begin{eqnarray*} 
&\mathcal{J}_{O,3,10}=\langle u,v \rangle,\ 
 \mathcal{J}_{O,4,10}=\langle u^3,v \rangle,\ 
 \mathcal{J}_{O,5,10}=\langle u^5,uv, v^2 \rangle,\ 
 \mathcal{J}_{O,6,10}=\langle u^7,u^3v, v^2 \rangle,\\ 
&\mathcal{J}_{O,7,10}=\langle u^{10},u^5v,uv^2,v^3 \rangle,\ 
 \mathcal{J}_{O,8,10}=\langle u^{12},u^7v,u^3v^2,v^3, h_2^{(2,0)}
 \rangle,\\ 
&\mathcal{J}_{O,9,10}=\langle u^{14},u^{10}v,u^5v^2,uv^3,v^4,r_2^{(4,0)}
 \rangle
\end{eqnarray*}
where  $h_2^{(2,0)}(u,v):=u^2h_2(u,v)$ and $r_2^{(4,0)}(u,v):=u^4r_2(u,v)$. 
\end{assertion}
\noindent
The proof  follows  from Lemma 3 and Lemma 4   and by  an easy
computation.  
 
Thus we have $\rho_8(O)=21$, $\rho_9(O)=29$ and 
\begin{equation*}
\tilde{\rho}(k)=
  \begin{cases}
    1&  k=7,9,\\
    0&  k\ne 7,9.
  \end{cases},\quad 
\iota_{k}(O)>10(k-3), \quad 3\le k \le 9.
\end{equation*}

\begin{assertion}
The map $\bar{\sigma}_k$ is injective for all $k=3,\dots,9$.
\end{assertion} 
\begin{proof}
 Recall that $C$ can be either irreducible or reducible in this case.
 As $\iota_{k}(O)>10(k-3)$, if  
 $C$ is irreducible, then the assertion follows from Proposition 1. 

\noindent
Assume $C$ is not irreducible. 
We have seen in the previous argument in \S3.4,
$C$ has two irreducible components  of respective degree 1 and 9.
Namely we can write $C=C_1 \cup C_9$ where $C_1=\{y=0\}$.
Suppose that  there exists a non-zero 
$g\in {\rm{Ker\,}}\bar{\sigma}_k\subset O(k-3)$.
As $\iota_{k}(O)>10(k-3)$, $g$ divides $f$
by the proof of Proposition 1.  
This is possible  only if $k\ge  4$ and  $\deg\,g=1$. 
By the assumption, we have $g=c\,y$
with $c\ne 0$. As $y=v+\varphi(u)$, we see that $g$ can not be in the ideal 
$\mathcal J_{O,k,10}$ for $k\ge  5$, as $v\notin \mathcal J_{O,k,10}$ by Assertion 1. This implies that
 $\bar{\sigma}_k$ is injective for $k\ne 4$.
Assume $k=4$.
As $a_{20}\ne 0$, ${\rm ord}_u\,\varphi(u)=2$ and 
$\mathcal{J}_{O,4,10}=\langle u^3,v \rangle$, again we see that 
 $v+\varphi(u)\notin \mathcal{J}_{O,4,10}$. 
This is a contradiction for $g\in {\rm{Ker\,}}\,\bar{\sigma}_4$
and  the proof is completed.
\end{proof}

Therefore we obtain the property  ($\sharp$):
  $\ell_k=1$ for $k=7,9$ and $\ell_k=0$ otherwise. Thus the
 reduced Alexander polynomial is given by  
 $\tilde{\Delta}_{C}(t)=t^4-t^3+t^2-t+1$
for the case $(C,O)\sim (B_{9,2}^2)^{B_{5,2}}\circ B_{2,1}$. 

\subsection{The case of $B_{29,2}\circ B_{2,1}\circ (B_{2,1}^2)^{B_{5,2}}$.}
By local classification \cite{Kawashima1}, 
this singularity 
appears in the case that the associated curves
$C_2$ and $C_5$  satisfies following  conditions:
\begin{enumerate}
\item $C_2$ is irreducible and $I(C_2,C_5;O)=10$. 

\item The multiplicity of $(C_5,O)$ is 3 and 
the tangent cone of $C_5$ at $O$
 consists of a line $L_1$ with multiplicity 2 and
             a single line $L_2$.

\item The conic $C_2$ is tangent to the line  $L_1$  at $O$. 
\end{enumerate}
Suppose that $C_2$ and $C_5$ satisfies the above conditions.
Then we may assume that the defining polynomials of $C_2$ and $C_5$ are 
following forms:
\begin{eqnarray*}
\begin{split}
f_2(x,y)=&y+a_{20}x^2+a_{11}xy+a_{02}y^2,\ a_{20}\ne 0\\
f_5(x,y)=&b_{05}y^5
         +a_{02}^2b_{12}\,xy^4
         +2a_{02}b_{12}x(a_{11}x+1)y^3\\&
         +\left(\frac{1}{27}b_{12}
(4a_{02}b_{12}^2+54a_{02}a_{20}+27a_{11}^2)x^3+2a_{11}b_{12}x^2+b_{12}x\right) {y}^{2}\\
&+\left(\dfrac{2}{27}b_{12}x^3(2b_{12}^2+27a_{20})(a_{11}x+1)\right)y
+\frac{1}{27}b_{12}a_{20}(27a_{20}+4b_{12}^2)x^5
\end{split}
\end{eqnarray*}
where $b_{12}\ne0$ and $b_{12}^2+9a_{20}\ne 0$ in general.
If $b_{12}^2+9a_{20}=0$, 
$C$ has the line component which is defined by $\{y=0\}$.
Now we take a local coordinates $(u,v)$ of the following type so that  
\begin{eqnarray*}
\begin{split}
 & x= u,\quad 
   y= v+\varphi(u),\quad \varphi(u)=-a_{20}u^2+\cdots,\ a_{20}\ne0,\\
 &    f_2(u, v+\varphi(u))=v+\psi(u)=v+\beta_7u^7+{\text{(higher terms)}}, \quad \beta_7\ne0, \\
 & 
 f_5(u, v+\varphi(u))=b_{05}\,v^5+b_{12}uv(v+\dfrac{4}{27}b_{12}^2u^2)+c_4\,u^{18}+{\text{(higher
 terms)}},\ b_{12}\ne0,\\ 
 & 
 f(u, v+\varphi(u))=v^2(v+d_1u^2)(v+d_2u^2)^2+\beta_7^5u^{35}+{\text{(higher
 terms)}},\ d_1,d_2\ne0.
 \end{split}
\end{eqnarray*}
By an explicit calculation, we have
$ d_2=\dfrac{4}{9}\,b_{12}^2$ and $d_2+a_{20}\ne 0$.
(If $d_2+a_{20}= 0$, $f$ becomes a non-reduced polynomial.)
Then the Newton boundary $\Gamma(f;u,v)$ consists of two faces 
$\Delta_1$ and  $\Delta_2$ so that 
their face functions are given by
\[
f_{\Delta_1}(u,v)=v^2(v+d_1u^2)(v+d_2u^2)^2,\quad 
f_{\Delta_2}(u,v)=u^6(d_1d_2^2v^2+\beta_7^5u^{29}).
\]
Note that $f(u,v)$ is degenerate on $\Delta_1$.
We take the  canonical toric modification  $\pi_1:X_1\to \Bbb{C}^2$
with respect to  $\{\Sigma_1^{*},(u,v),O\}$ where
$\Sigma_1^{*}$ is the canonical regular simplicial cone subdivision 
with vertices
\begin{equation*}
E_1,\quad
Q_k = \begin{pmatrix}
         1  \\
         k 
         \end{pmatrix}(1\le k \le 14),\quad
Q_{15} = \begin{pmatrix}
         2  \\
         29 
         \end{pmatrix},\quad
Q_{16} = \begin{pmatrix}
         1  \\
         15 
         \end{pmatrix},  \quad
E_2              
\end{equation*}
where $Q_2$ and $Q_{15}$ 
are  the weight vectors of the faces $\Delta_1$ and $\Delta_2$ respectively. 
Then  the divisor $(\pi_1^{*}f)$ is given by 
\[(\pi_1^{*}f)= \tilde{C}+\sum_{i=1}^{16}m(f,Q_i)\hat{E}(Q_i),\]
where $\tilde{C}$ is the strict transform of $C$ 
and  intersects only with the exceptional divisors $\hat{E}(Q_2)$ and
$\hat{E}(Q_{15})$.
We can see that $\tilde{C}$ is smooth  at $\tilde{C}\cap
\hat{E}(Q_{15})$
and the intersection is transverse.
 On the other hand,
$\tilde{C}$ intersects with $\hat{E}(Q_{2})$
at two points $\xi_{1,1},\, \xi_{1,2}$ 
where  $\xi_{1,1}=(0,-d_1),\, \xi_{1,2}=(0,-d_2)$
in the toric coordinates  $(u_1,v_1)$ of the chart $\Bbb{C}_{\tau}^2$ with
$\tau={\rm{Cone}} \,(Q_2,Q_3)$.
Note that 
$(\tilde{C},\xi_{1,1})$ is smooth and the intersection with
$\hat{E}(Q_{2})$ 
 is transverse
at  $\xi_{1,1}=(0,-d_1)$.
On the other hand,
$(\tilde{C},\xi_{1,2})$ has singularity. 
To see the singularity $(\tilde{C},\xi_{1,2})$,
we take the admissible translated coordinates $(u_1,v_1')$ with
$v_1'=v_1+d_2+h(u_1)$ where $h$  takes the form
$h(u_1)=q_1u+q_2u^2$.
Then we see that 
$\pi_1^*f(u_1,v_1')=c\,u_1^{10}({v_1'}^2+\beta u_1^5+{\text{(higher
terms)}})$
and $(\tilde{C},\xi_{1,2})\sim B_{5,2}$.
Now we take the second toric modification $\pi_2:X_2\to X_1$
with respect to $\{\Sigma_2^{*},(u_1,v'_1),\xi_{1,2}\}$ where 
 $\Sigma_2^{*}$ is  the canonical regular 
simplicial cone subdivision with vertices
\begin{equation*}
E_1,\quad 
T_1 = \begin{pmatrix}
         1  \\
         1 
         \end{pmatrix},\
T_2 = \begin{pmatrix}
         1  \\
         2 
         \end{pmatrix},\
T_3 = \begin{pmatrix}
         2  \\
         5 
         \end{pmatrix},  \
T_4 = \begin{pmatrix}
         1  \\
         3
         \end{pmatrix},  \quad
E_2
\end{equation*}
where the weight vector $T_3$ 
corresponds to the unique face of $\Gamma(\pi_1^{*}f;u_1,v_1')$.
Then we  have
\begin{multline*}    
(\pi^{*}\,f\,)=5\hat{E}(Q_1)+\sum_{i=2}^{14}2(i+3)\hat{E}(Q_i)
             +70\hat{E}(Q_{15})+35\hat{E}(Q_{16})\\
            \hspace{7cm}+12\hat{E}(T_1)+14\hat{E}(T_2)+30\hat{E}(T_3)
              +15\hat{E}(T_4).\\
\quad (\pi^{*}K)=\hat{E}(Q_1)+\sum_{i=2}^{14}i\hat{E}(Q_i)
           +30\hat{E}(Q_{15})+15\hat{E}(Q_{16})
           +3\hat{E}(T_1)+4\hat{E}(T_2)+10\hat{E}(T_3)
           +5\hat{E}(T_4)
\end{multline*}
and we consider two polynomials $h_1(u,v)$ and $r_1(u,v)$
  which are defined by
 $h_1(u,v)=v+d_2u^2$ and $r_1(u,v)=h_1(u,v)-\frac{q_1}{d_2}u^3$.
  Then
  \begin{eqnarray*}
   & \pi_1^*h_1(u_1,v_1')=u^{2}(e_3v'_1+e_4u_1+\text{(higher terms)}),\\
   & \pi_1^*r_1(u_1,v_1')=u^{2}(e'_3v'_1+e'_4u_1^2+\text{(higher terms)}),\\
   & m(h_1,Q_2)=2,\ \ m(h_1,Q_{15})=4,\ \ m(h_1,T_3)=6,\\
   & m(r_1,Q_2)=2,\ \ m(r_1,Q_{15})=4,\ \ m(r_1,T_3)=8.
  \end{eqnarray*}

\begin{remark}{\rm{The Alexander polynomial does not changes on the
 irreducible  component of the configuration space of the fixed 
topological type of the singularity.
Therefore 
for the practical computation, it is easier to choose some explicit
 values.
We take $b_{05}=b_{12}=a_{02}=a_{20}=a_{11}=1$. 
Then we have
\begin{eqnarray*}
\begin{split}
f_2(x,y)&=y^2+(x+1)y+x^2\\
f_5(x,y)&=y^5+xy^4+ 2(x^2+x)y^3+(\frac {85}{27}x^3+2x^2+x)y^2
         +(\frac{58}{27}x^4+\frac {58}{27}x^3)y+\frac {31}{27}x^5\\
\varphi(u)&=-u^2+u^3-2u^4+4u^5-9u^6+\frac {111}{4}u^7-\frac{183}{2}u^8
         +316u^9-1079u^{10}+\frac{7259}{2}u^{11}\\
 &\qquad\qquad\qquad\qquad\qquad -\frac{801559}{64}u^{12}+\frac{2872109}{64}u^{13}
      -\frac {5348333}{32}u^{14}.
\end{split}
\end{eqnarray*}
}}\end{remark}
\begin{assertion}
Under the above situation,
\begin{enumerate}
\item[(a)] 
The ideals $\mathcal{J}_{O,k,10}$ are given by 
\begin{eqnarray*} 
& \mathcal{J}_{O,3,10}=\langle u,v \rangle,\quad 
  \mathcal{J}_{O,4,10}=\langle u^2,v \rangle,\quad 
  \mathcal{J}_{O,5,10}=\langle u^3,uv, v^2 \rangle,\quad 
  \mathcal{J}_{O,6,10}=\langle u^6,u^2v, v^2 \rangle,\\
& \mathcal{J}_{O,7,10}=\langle u^{10},u^4v,u^2v^2,v^3,
                       h_1^{(1,1)} \rangle,\quad 
  \mathcal{J}_{O,8,10}=\langle u^{13},u^5v,u^3v^2,uv^3,v^4, 
                               h_1^{(2,1)},
                               h_1^{(0,2)}
                       \rangle,\\ 
& \mathcal{J}_{O,9,10}=\langle u^{17}, u^7v, u^5v^2, u^3v^3, uv^4, v^5, 
                               r_1^{(3,1)}, 
                               r_1^{(1,2)}, 
                               h_1^{(0,3)}
                       \rangle,
\end{eqnarray*} 
where $h_1^{(r,s)}(u,v):= u^rv^s h_1(u,v)$ and
      $r_1^{(r,s)}(u,v):= u^rv^s r_1(u,v)$.
\item[(b)] 
The kernel of $\bar\sigma_{k}$ are given by
\begin{eqnarray*}
& {\rm{Ker\,}}\bar\sigma_{3} = \langle 0 \rangle,\quad 
  {\rm{Ker\,}}\bar\sigma_{4} = \langle y \rangle,\quad
  {\rm{Ker\,}}\bar\sigma_{5} = \langle y^2,xy \rangle,\quad 
  {\rm{Ker\,}}\bar\sigma_{6} = \langle yf_2,y^3 \rangle, \\
& {\rm{Ker\,}}\bar\sigma_{7} = \langle y^2f_2 \rangle,\quad 
  {\rm{Ker\,}}\bar\sigma_{8} = \langle y^3f_2 \rangle,\quad 
  {\rm{Ker\,}}\bar\sigma_{9} = \langle 3yf_5-2b_{12}xyf_2^2-c_0y^2f_2^2
                               \rangle.\\
&\qquad \text{with} \qquad c_0=\dfrac{81a_{11}a_{20}}{4b_{12}(9a_{20}+4b_{12}^2)}.
\end{eqnarray*}
\end{enumerate} 
\end{assertion}

\noindent {\em Proof.}
The assertion (a) follows from Lemma 3 and Lemma 4.  
We consider the assertion (b). 
By the choice of the local coordinates $(u,v)$, 
we have relations:
\begin{equation*}
\begin{split}
 & x= u, \quad
  y= v+\varphi(u)=v-a_{20}u^2+\cdots, \\
&  f_2(u,v+\varphi(u))=v+\psi(u)=v+\beta_7u^7+{\text{(higher terms)}},\ \beta_7\ne 0,\\
&  f_5(u,v+\varphi(u))=b_{05}\,v^5
        +b_{12}uv(v+\dfrac{4}{27}b_{12}^2u^2)+c_4\,u^{18}+{\text{(higher terms)}}.
\end{split}
\end{equation*}
Put $\varphi(u)=\sum_{j=2}^{\infty}\alpha_j u^j$ with $\alpha_2=-a_{20}$.
We define
\[
\ord\, \mathcal{J}_{O,k,10}:=\min\{\ord_{(u,v)}\, h \,|\,h\in 
                                         \mathcal{J}_{O,k,10}\}.
\]
Thus for any $g\in {\rm{Ker\,}}\bar\sigma_{k}$, 
we have 
\begin{eqnarray}
\label{order-estimate}\quad 
 \ord_{(x,y)}\, g =\ord_{(u,v)}\sigma_k(g)\ge 
 \ord\, \mathcal{J}_{O,k,10}.
\end{eqnarray}

\noindent{\bf{Case $k=4$}}: ${\rm{Ker\,}}\bar\sigma_{4}=\langle y
 \rangle$.
\begin{proof}
The inclusion 
$\langle y \rangle\subset {\rm{Ker\,}}\bar\sigma_{4}$
 holds by the definition of $\sigma_{4}$.
For any $g\in {\rm{Ker\,}}\bar\sigma_{4}\subset O(1)$,
writing 
$g(x,y)=c_1+c_2x+c_3y$,
\[
 \sigma_4(g)(u,v)= c_1+c_2u+c_3(v+\varphi(u))\in
 \mathcal{J}_{O,4,10}=\langle u^2,v \rangle.
\]
Hence we have $c_1=c_2=0$ and 
${\rm{Ker\,}}\bar\sigma_{4}\subset\langle y \rangle$.
\end{proof}

\noindent{\bf{Case $k=5$}}:
${\rm{Ker\,}}\bar\sigma_{5}=\langle y^2,xy \rangle$.
\begin{proof}
First we show that $y^2,\,xy\in {\rm{Ker\,}}\bar\sigma_{5}$.
By the definition of $\sigma_{5}$, we have 
\begin{equation*}
\begin{split}
\sigma_5(y^2)=&(v+\varphi(u))^2 
             =v^2-2a_{20}u^2v+a_{20}^2u^4+
              {\text{(higher terms)}}
              \in \mathcal{J}_{O,5,10}\\
\sigma_5(xy)=&u(v+\varphi(u))
             =uv-a_{20}u^3+
             {\text{(higher terms)}}\in \mathcal{J}_{O,5,10}
\end{split}
\end{equation*}
as $\mathcal{J}_{O,5,10}=\langle u^3,uv, v^2 \rangle$. 
Next we show that 
${\rm{Ker\,}}\bar\sigma_{5}\subset \langle y^2,xy \rangle$.
Take $g\in {\rm{Ker\,}}\bar\sigma_{5}\subset O(2)$.
As $\ord\,\mathcal{J}_{O,5,10}=2$,
we can write 
$g(x,y)=c_1x^2+c_2xy+c_3y^2$
by $(\ref{order-estimate})$ and 
\[
\sigma_5(g)(u,v)=c_1u^2+c_2uv+c_3v^2+{\text{(higher terms)}}
\in \mathcal{J}_{O,5,10}=\langle u^3,uv,v^2 \rangle.
\]
Hence we have $c_1=0$ and
${\rm{Ker\,}}\bar\sigma_{5}=\langle y^2,xy \rangle$.
\end{proof}

\noindent{\bf{Case $k=6$}}:
${\rm{Ker\,}}\bar\sigma_{6}=\langle yf_2, y^3\rangle$.
\begin{proof}
First we show that $yf_2,\, y^3\in {\rm{Ker\,}}\bar\sigma_{6}$.
By the definition of $\sigma_{6}$, we have 
\begin{equation*}
\begin{split}
\sigma_6(yf_2) =& (v+\varphi(u))(v+\psi(u))\\
               =& v^2-a_{20}u^2v-a_{20}\beta_7u^9+
                  {\text{(higher terms)}}
                  \in \mathcal{J}_{O,6,10}\\
\sigma_6(y^3) =& (v+\varphi(u))^3
              =  v^3-3a_{20}u^2v^2+3a_{20}^2u^4v-a_{20}^3u^6
                  +{\text{(higher terms)}}
                  \in \mathcal{J}_{O,6,10}
\end{split}
\end{equation*}
as $\mathcal{J}_{O,6,10}=\langle u^6,u^2v, v^2 \rangle$.
 
Next we show that ${\rm{Ker\,}}\bar\sigma_{6}\subset \langle yf_2, y^3\rangle$.
Take  $g\in {\rm{Ker\,}}\bar\sigma_{6}\subset O(3)$.
As $\sigma_6(g)\in \mathcal{J}_{O,6,10}$, we can write
 \[
  \sigma_6(g)(u,v)=g(u,v+\varphi(u))=a_1(u,v)u^6+a_2(u,v)u^2v+a_3(u,v)v^2
 \]
where $a_i\in \mathcal{O}_O$ $(i=1,2,3)$.
Define $g'(u,v)$ by the above right side polynomial. 
Then we see that 
\[
 I(y,g;O)=\ord_u g'(u,-\varphi(u)) \ge 4.
\] 
On the other hand, if $y$ does not divide $g$, 
$I(y,g;O)\le 3$ by  B\'ezout's theorem which is an obvious contradiction.
Therefore  $y$ divides $g$. 
 Thus we can write $g(x,y)=yg_2(x,y)$ where
 $g_2\in O(2)$.
Dividing $g_2$ by $f_2$ as a polynomial of $x$, we can write $g_2$ as 
$g_2=c_0f_2+(c_1+c_2y)x+c_3y^2+c_4y+c_5$
for some constants $c_0,\dots, c_5$.
As $y f_2,\,y^3\in {\rm{Ker\,}}\bar\sigma_{6}$,
we need to have $y((c_1+c_2y)x+c_4y+c_5)\in {\rm{Ker\,}}\bar\sigma_{6}$.
By a simple computation, we conclude $c_1=c_2=c_4=c_5=0$ and 
\[
g(x,y)=c_0 yf_2(x,y)+c_3 y^3\in \langle yf_2, y^3\rangle.
\]
\end{proof}

\noindent{\bf{Case $k=7$}}: 
${\rm{Ker\,}}\bar\sigma_{7}=\langle y^2f_2\rangle$.
\begin{proof}
First we show that 
$y^2f_2\in {\rm{Ker\,}}\bar\sigma_{7}$.
By the definition of $\sigma_7$, we have 
\begin{equation*}
\begin{split}
\sigma_7(y^2f_2)(u,v)=& (v+\varphi(u))^2(v+\psi(u))\\
                     =& v^3-2a_{20}u^2v^2+a_{20}^2u^4v
                        +a_{20}^2\beta_7u^{11}
                        +{\text{(higher terms)}}
                        \in \bar{\mathcal{J}}_{O,7,10}
\end{split}
\end{equation*}
as
\begin{eqnarray*}
\begin{split}
& \bar{\mathcal{J}}_{O,7,10}=\langle u^{10},u^4v,u^2v^2,v^3\rangle,\quad
\mathcal{J}_{O,7,10}=\langle u^{10},u^4v,u^2v^2,v^3,
                     h_1^{(1,1)} \rangle
\end{split}
\end{eqnarray*}
where 
$h_1^{(1,1)}(u,v):= uv(v+d_2u^2)$.
Next we show that 
${\rm{Ker\,}}\bar\sigma_{7}\subset
            \langle y^2f_2\rangle$.              
Take $g\in {\rm{Ker\,}}\bar\sigma_{7}\subset O(4)$ and
we can write
$\sigma_7(g)$ as
\begin{eqnarray*}
 \begin{split}
& \sigma_7(g)(u,v)=\sum_{i\ge 0}g_i(u)v^i,\quad  \ord_u g_0(u)\ge 10, \ 
\ord_u g_1(u)\ge 3,\
  \ord_u g_2(u)\ge  1.
 \end{split}
\end{eqnarray*}
Then we see that  
$I(g,y;O)=\ord_u\,\sigma_7(g)(u,-\varphi(u))\ge 5$ and
by B\'ezout's theorem,  $y$ divides $g$.
Similarly we can see that 
we have $I(g,f_2;O)=\ord_u\,\sigma_7(g)(u,-\psi(u))\ge 10$ and
again by B\'ezout's theorem, 
we conclude $f_2$ divides $g$.  
Thus we can write $g(x,y)=yf_2(c_0+c_1x+c_2y)$
for some $c_0,c_1,c_2\in \Bbb C$.
The assumption $g,\,y^2f_2\in {\rm{Ker\,}}\bar\sigma_{7}$
implies that 
$g(x,y)-c_2\,y^2f_2=(c_0+c_1x)yf_2\in {\rm{Ker\,}}\bar\sigma_{7}$.
Thus we have 
$\sigma_7(c_0yf_2)(u,0)=-c_0a_{20}\beta_7u^9+\cdots\in
\mathcal{J}_{O,7,10}$.
Therefore 
$c_0=0$ as $\ord_u \sigma_7(c_0yf_2)(u,0)\ge 10$. 
Moreover we have
 \[
 \begin{split}
\sigma_7(g)\equiv
\sigma_7(c_1xyf_2)&\equiv c_1\,uv(v-a_{20}u^2) \mod
  \bar{\mathcal{J}}_{O,7,10}.
\end{split}
\]
As $d_2+a_{20}\ne 0$, we see that 
$uv(v-a_{20}u^2)\notin\mathcal{J}_{O,7,10} $.
Hence we have $c_1=0$ and
 we conclude $g(x,y)=c_2\,y^2f_2$.
\end{proof}

\noindent{\bf{Case $k=8$}}: ${\rm{Ker\,}}\bar\sigma_{8}=\langle y^3f_2\rangle$.
\begin{proof}
First we show that $y^3f_2\in {\rm{Ker\,}}\bar\sigma_{8}$.
By the definition of $\sigma_8$, we have 
\begin{equation*}
\begin{split}
\sigma_8(y^3f_2)(u,v)= &(v+\varphi(u))^3(v+\psi(u))\\
                     = &v^4-3a_{20}u^2v^3+3a_{20}^2u^4v^2
                       -a_{20}^3u^6v-a_{20}^3\beta_7u^{13}
                       +{\text{(higher terms)}}.
\end{split}
\end{equation*}
As 
$\bar{\mathcal{J}}_{O,8,10}=\langle u^{13},u^5v,u^3v^2,uv^3,v^4\rangle$,
we see that $\sigma_8(y^3f_2)  \in \bar{\mathcal{J}}_{O,8,10}$. 

Next we show that 
${\rm{Ker\,}}\bar\sigma_{8}
\subset\langle y^3f_2\rangle$.              
Take  $g\in {\rm{Ker\,}}\bar\sigma_{8}\subset O(5)$ and
 write
$\sigma_8(g)$ as
\begin{eqnarray*}
 \begin{split}
& \sigma_8(g)(u,v)=\sum_{i\ge 0}g_i(u)v^i,\quad  \ord_ug_0(u)\ge 13,\  
  \ord_ug_1(u)\ge 4, \ 
  \ord_ug_2(u)\ge 2.
\end{split}
\end{eqnarray*}
As $I(g,y;O)=\ord_u\,\sigma_8(g)(u,-\varphi(u))\ge 6$,
 we see that $y$ divides $g$
 by B\'ezout's theorem.
Similarly we can see that 
   $I(g,f_2;O)=\ord_u\,\sigma_8(g)(u,-\psi(u))\ge 11$,  
 we see that  $f_2$ divides $g$.
Hence we have $g(x,y)=yf_2g'(x,y)$ for some $g'\in O(2)$. 
We put $g'(x,y)=c_{02}y^2+r(x,y)$ where
 $r(x,y)=c_{11}xy+c_{01}y+c_{20}x^2+c_{10}x+c_{00}$.
As $y^3f_2\in {\rm{Ker\,}}\bar\sigma_{8} $, we have
$yf_2r(x,y)\in {\rm{Ker\,}}\bar\sigma_{8}$.
Consider the expression 
\begin{eqnarray*}
 \begin{split}
& \sigma_8(yf_2r)(u,v)=\sum_{i\ge 0}\psi_i(u)v^i
\end{split}
\end{eqnarray*}
We can see that 
$\psi_2(u)=c_{00}+(c_{00}a_{20}-c_{10})u+u^2\tilde{\psi}_2(u)$.
Thus $c_{00}=0$ and $c_{10}=0$  as $\ord_u\psi_2(u)\ge 2$. 
Now  we have
 \[
\psi_0(u)=-a_{20}\beta_7(a_{20}c_{01}-c_{20})u^{11}
                   +a_{20}^2
 \left( d(a_{20}c_{01}-c_{20})+\beta_7(a_{11}c_{01}-c_{11})\right)u^{12}
                   +u^{13}\tilde{\varphi}_0(u).
\]
As $\ord_u \psi_0(u)\ge 13$ by the assumption $\sigma_8(yf_2r)\in
 \mathcal{J}_{O,8,10}$, 
the coefficients of $u^{11},u^{12}$ in $\psi_0(u)$ must vanish.
Therefore
 $a_{20}c_{01}-c_{20}=0$ and
 $a_{11}c_{01}-c_{11}=0$. 
Thus we conclude $r(x,y)=c_{01}(y+a_{20}x^2+a_{11}xy)$.
Consider the weight vector $P={}^t(1,2)$ for the variables $u,v$.
Then we compute the leading term of $\sigma_{8}(yf_2r)(u,v)$
with respect to $P$:
\[
 \sigma_{8}(yf_2r)_P(u,v)=
          c_{01}v^2(v-a_{20}u^2)
\]
As the lowest degree of elements in  $\mathcal{J}_{O,8,10}$
is $6$ and they are generated by $h_1^{(2,1)}=u^2v(v+d_2u^2)$ and 
$h_1^{(0,2)}=v^2(v+d_2u^2)$. Thus we must have the equality
$\sigma_{8}(yf_2r)_P(u,v)|_{v=-d_2u^2}=0$. This implies that $c_{01}=0$
as $d_2+a_{20}\ne 0$.
\end{proof}

\noindent{\bf{Case $k=9$}}: 
${\rm{Ker\,}}\bar\sigma_{9}= 
\langle
  3yf_5-2b_{12}xyf_2^2-c_0y^2f_2^2
 \rangle$
where $c_0=\dfrac{81a_{11}a_{20}}{4b_{12}(9a_{20}+4b_{12}^2)}$.
\begin{proof}
The proof of this case is most computational. As the Alexander polynomial is a
 topological invariant,  we can choose any polynomial in the connected
component of the moduli space. Thus we use 
 the  polynomial in Remark 2.
(We take $b_{05}=b_{12}=a_{02}=a_{20}=a_{11}=1$.)
\begin{eqnarray*}
\begin{split}
f_2(x,y)&=y^2+(x+1)y+x^2\\
f_5(x,y)&=y^5+xy^4+ 2(x^2+x)y^3+\left(\frac {85}{27}x^3+2x^2+x\right)y^2
         +\left(\frac{58}{27}x^4+\frac {58}{27}x^3\right)y+\frac {31}{27}x^5\\
\varphi(u)&=-u^2+u^3-2u^4+4u^5-9u^6+\frac {111}{4}u^7-\frac{183}{2}u^8
         +316u^9-1079u^{10}+\frac{7259}{2}u^{11}\\
 &\qquad\qquad\qquad\qquad\qquad -\frac{801559}{64}u^{12}+\frac{2872109}{64}u^{13}
      -\frac {5348333}{32}u^{14}.
\end{split}
\end{eqnarray*}
and then $h_1(u,v)=v+\frac{4}{9}u^2$,
$r_1(u,v)=
 (v+\frac{4}{9}u^2)-\frac{13}{9}u^3$ and
$c_{0}=\frac{81}{52}$.

First we show that
$p(x,y):=3yf_5-2xyf_2^2-\frac{81}{52}y^2f_2^2$
is in the kernel of $\bar\sigma_{9}$.
Put $p_1(x,y)=3yf_5-2xyf_2^2$.
We  observe that  
\begin{equation*}
\begin{split}
\sigma_9(p_1)(u,v)         &\equiv uv(v-u^2)(v+\dfrac{4}{9}u^2)+
                             u^2v^2\left(2v-\dfrac{5}{9}u^2\right)
                         \mod \bar{\mathcal{J}}_{O,9,10}.\\
\sigma_9(\dfrac{81}{52}y^2f_2^2) &\equiv \dfrac{81}{52}
 v^2(v-u^2)^2
                         \mod \bar{\mathcal{J}}_{O,9,10}.\\
\sigma_9(p)(u,v)      
   &\equiv -\left(1+\dfrac{13}{4}u\right)r_1^{(3,1)}
  +\left(1+\dfrac{151}{26}u\right)r_1^{(1,2)}
  -\dfrac{81}{52}h_1^{(0,3)}
  \ \ \quad  \mod \bar{\mathcal{J}}_{O,9,10}.
\end{split}
\end{equation*}    
Thus we conclude $\sigma_9(p)\in {\mathcal{J}}_{O,9,10}$, 
as $\mathcal{J}_{O,9,10}=
\langle  
        u^{17}, u^7v, u^5v^2, u^3v^3, uv^4, v^5, 
        r_1^{(3,1)}, 
        r_1^{(1,2)}, 
        h_1^{(0,3)}
\rangle$.  

Next we will show that 
${\rm{Ker\,}}\bar\sigma_{9}$ is generated by $p$.
Take  $g \in {\rm{Ker\,}}\bar\sigma_{9}\subset O(6)$ with $g\ne 0$.
As $\ord\, \mathcal{J}_{O,9,10}=4$, we have $\ord_{(x,y)}\, g=4$.
Hence we can put 
\[
 g(x,y)=\sum_{4\le r+s \le 6}c_{rs}x^ry^s,\quad
 \sigma_9(g)(u,v)=\sum_{i\ge 0}\psi_i(u)v^i,
\] 
where 
\[
\psi_0(u)=\sum_{i\ge 4}a_i u^i,\quad 
\psi_1(u)=\sum_{i\ge 3}b_i u^i,\ \
{\text{$a_i$,  $b_i$ are linear polynomials in $c_{rs}$.}} 
\]
By the assumption 
$\sigma_9(g)\in \mathcal{J}_{O,9,10}$,
we have $\ord_u\psi_0(u)\ge 17$ and $\ord_u\psi_1(u)\ge 5$.
Hence we solve the 15-equations 
$a_i=0$, $4\le i\le 16$ and $b_i=0$, $i=3,4$ 
in $c_{rs}$ for the lexicographical order.
After solving these equations, $g$ takes the form:
\begin{equation*}
\begin{split}
g(x,y)=&c_{06}y^6+(c_{15}x+c_{05})y^5
      +\left((2c_{15}-\dfrac{1}{2}c_{05})x^2+(2c_{15}-c_{05})x
      +\dfrac{1}{2}c_{05}\right)y^4\\
    &+\Big{(}(\dfrac{4}{27}c_{06}-\dfrac{56}{27}c_{05}+3c_{15})x^3+
         (2c_{15}-c_{05})x^2+(c_{15}-c_{05})x\Big{)}y^3\\
    &+\Big{(}(\dfrac{4}{27}c_{06}-\dfrac{85}{54}c_{05}+2c_{15})x^4
     +(2c_{15}+\dfrac{4}{27}c_{06}-\dfrac{56}{27}c_{05})x^3\Big{)}y^2\\&+
     \left(c_{15}-\dfrac{29}{27}c_{05}+\dfrac{4}{27}c_{06}\right)x^5y
\end{split}
\end{equation*}
where $g$ has still 3-parameters $c_{05},\,c_{15},\, c_{06}$.
We get the equality 
\begin{equation*}
\begin{split}
\sigma_{9}(g)(u,v)\equiv &
(c_{15}-c_{05})uv^3
+(\dfrac{4}{27}c_{06}-c_{15}+\dfrac{25}{27}c_{05})u^3v^2
+(\dfrac{2}{27}c_{05}-\dfrac{4}{27}c_{06})u^5v \\
& +\dfrac{1}{2}c_{05}v^4+(2c_{15}-3c_{05})u^2v^3
+(\dfrac{4}{27} c_{06}-c_{15}+\dfrac{77}{54}c_{05})u^4v^2
\mod \bar{\mathcal{J}}_{O,9,10}.
\end{split}
\end{equation*}
Now we consider the weight vector $P={}^t(1,2)$
 for variables $u,v$ as in the previous case.
Then $\deg_P\sigma_9(g)(u,v)=7$ and 
 \[
 \begin{split}
\sigma_9(g)_P (u,v)&=
(c_{15}-c_{05})uv^3
+(\dfrac{4}{27}c_{06}-c_{15}+\dfrac{25}{27}c_{05})u^3v^2
+(\dfrac{2}{27}c_{05}-\dfrac{4}{27}c_{06})u^5v,\\
 \end{split}
 \]
As the lowest degree of the generators of $\mathcal{J}_{O,9,10}$ is also
 $7$
and they are $r_1^{(3,1)}$ and $r_1^{(1,2)}$.
Thus we must have
\[
 \sigma_9(g)_P (u,v)=( a r_1^{(3,1)}+b r_1^{(1,2)})_P,\quad \text{for some}\,\,
 a,b\in \Bbb C.
\]
Thus we see that $\sigma_9(g)_P(u,-\frac{4}{9} u^2)=0$.
This gives the equality 
 \[
 5c_{05}-6c_{15}+2c_{06}=0
\]
and we eliminate  $c_{05}$ using the above equality and then
$a,b$ are solved as follows.
\[
  a= \frac{1}{5}\,{ c_{15}}-\frac{2}{5}\,{ c_{06}},\quad 
  b=-\frac 15\,{c_{15}}+\frac 25\,{ c_{06}}.
\]
We put $g_1=\sigma_9(g)-(ar_1^{(3,1)}+br_1^{(1,2)})$.
Then we see that $\deg_Pg_1=8$.
Thus we can write 
\[
g_1(u,v)=a'ur_1^{(3,1)}+b'ur_1^{(1,2)}+c'h_1^{(0,3)},\quad 
{\text{for some}}\ \  a',b',c'\in \Bbb{C}.		   
		  \]
Again we need to have $g_1(u,-4/9u^2)=0$ which gives the equality
\[
 15c_{15}+22c_{06}=0.
\]
Eliminating the parameter $c_{15}$, we  finally obtain the expression
\begin{equation*}
\begin{split}
g(x,y)&=\dfrac{c_{06}}{675}y
 \Big{(}675 y^5-(990x+1458) y^4-(1251x^2+522x+729) y^3
 \\&\hspace{3cm}
+(154x^2-522x+468)x y^2+(415x^4+1144x^3) y+676 x^5\Big{)}
\end{split}
\end{equation*}
and  we conclude that
\[
g(x,y)=\dfrac{52}{75}c_{06}\,p(x,y).
\]
\end{proof}
The proof of Assertion 3 is  now completed.

\vspace{0.5cm}
Now we are ready to compute the Alexander polynomial for the case
$(C,O)\sim\linebreak
 B_{29,2}\circ B_{2,1}\circ (B_{2,1}^2)^{B_{5,2}}$.
By above assertions, we have 
$\rho_7(O)=15$, $\rho_8(O)=20$ and $\rho_9(O)=28$ hence 
 we obtain the property $(\sharp)$:
$\ell_k=1$ for $k=7,9$ and $\ell_k=0$ otherwise.
Therefore we  have $\tilde{\Delta}_{C}(t)=t^4-t^3+t^2-t+1$
by  Lemma 1. 
Thus the proof of Theorem 1 is completed.\\

\subsection{Linear torus curve}
The singularity $B_{50,2}$   appears also as a  linear torus curve of
type (5,2):
\[
 C:\quad f_5(x,y)^2-y^{10}=0
\]
with $I(f_5,y;O)=5$ (\cite{BenoitTu}).
In this case, $C$ consists of two smooth quintics and 
the Alexander polynomial  is given by following (\cite{BenoitTu}): 
\[
\Delta_{C}(t)=\frac{(t^{10}-1)}{t+1}.
\]
\subsection{Proofs of Corollary 1 and  Corollary 2}
The assertion of Corollary 1 is an immediate consequence of 
the Sandwich principle.  
The assertion of Corollary 2 is a result of \cite{BenoitTu}. In fact,
we only need to observe that 
the equivalence class of such  torus curves correspond 
bijectively to the partitions of 10 by locally intersection numbers 
of $C_2$ and $C_5$.
In particular, such a curve degenerates into an irreducible torus curve
with a unique singularity $B_{50,2}$ which corresponds to the partition 
$10=10$.

\bibliographystyle{abbrv}
\def\cprime{$'$} \def\cprime{$'$} \def\cprime{$'$} \def\cprime{$'$}
  \def\cprime{$'$}

\end{document}